%% file: 0paper.tex
\documentclass[12pt]{article}
\usepackage{amsfonts}
\usepackage{epsfig}

\title{The Triangular Bi-Pyramid Minimizes
a Range of Power Law Potentials}
\author{Richard Evan Schwartz \thanks{\hskip 5 pt Supported by 
N.S.F. Research Grant DMS-1204471}}

\newtheorem{theorem}{Theorem}[section]

\newtheorem{lemma}[theorem]{Lemma}

\newtheorem{corollary}[theorem]{Corollary}

\def\startproof{{\bf {\medskip}{\noindent}Proof: }}

\def\endproof{$\spadesuit$  \newline}

\def\D{\mbox{\boldmath{$D$}}}%
\def\N{\mbox{\boldmath{$N$}}}%
\def\Q{\mbox{\boldmath{$Q$}}}%
\def\R{\mbox{\boldmath{$R$}}}%
\def\Z{\mbox{\boldmath{$Z$}}}%

\begin{document}
\maketitle

\begin{abstract}
Combining a brilliant observation
of A. Tumanov with a divide-and-conquer
approach to Thomson's $5$-electron
problem, we give a rigorous computer-assisted
proof that the triangular bi-pyramid
is the unique
minimizer with respect to the potential
$$R_s(r)={\rm sign\/}(s)\frac{1}{r^s}$$
amongst all spherical $5$ point configurations
when  $s \in (-2,13]-\{0\}$.   The case
$s=-1$ corresponds to Polya's problem
and the case $s=1$ corresponds to
Thomson's problem.  The lower bound of
$13$ is fairly close to the presumed
optimal cutoff of $15.040809...$.
\end{abstract}

\input{1intro}

\input{2config}
\input{3sphere}

\input{4energy}
\input{5divide}
\input{6hessian}

\input{7poly}
\input{8tumanov}
\input{9compute}
\input{refs}

\end{document}

%% file: 1intro.tex
\section{Introduction}

\subsection{The Energy Minimization Problem}

Let $S^2$ denote the unit sphere in $\R^3$ and
let $P=\{p_1,...,p_n\}$ be a finite list of pairwise
distinct points on $S^2$.  
Given some function $f: (0,2] \to \R$ we can
form the energy potential
\begin{equation}
{\cal E\/}_f(P)=\sum_{i<j} f(\|p_i-p_j\|).
\end{equation}
For fixed $f$ and $n$,
one can ask which configuration(s)
minimize ${\cal E\/}_f(P)$. 

For this problem, the energy functional $f=R_s$, where
\begin{equation}
R_s(r)={\rm sign\/}(s)\frac{1}{r^s},
\end{equation}
is a natural one to consider.
At least for $s>0$, this is called
the {\it Riesz energy\/}, and we will
also use this name when $s<0$.
The case $s=1$, where $f(r)=1/r$,
 corresponds to the
electrostatic potential. This
case is known as {\it Thomson's problem\/}.
See [{\bf Th\/}].
The case $s=-1$ corresponds
to the problem of placing points
on the sphere so as to maximize
the total sum of the distances
between pairs. This is known as
{\it Polya's problem\/}. 

There is a large literature on the energy 
minimization problem. See [{\bf C\/}] for some early local results.
The online website [{\bf CCD\/}] is a compilation of experimental
results.
The paper [{\bf SK\/}] gives a nice survey in the two dimensional case,
with an emphasis on the case when $n$ is large. 
See also [{\bf RSZ\/}].
The paper [{\bf BBCGKS\/}] gives a survey
of results, both theoretical and experimental, about
highly symmetric configurations in higher dimensions.

When $n=2,3,4,6,12$, the most symmetric
configurations are the unique 
minimizers for all $R_s$
with $s \in (-2,\infty)-\{0\}$.  For the cases
$n=4,6$ (tetrahedron, octahedron) see [{\bf Y\/}].  For the case $n=12$
(icosahedron) see [{\bf A\/}].  All these cases are
covered by the vast result in [{\bf CK\/}, Theorem 1.2].

The case $n=5$ has been notoriously intractable.
In this case, there is a general feeling that 
for a wide range of energy choices, and in particular
for the Riesz energies, that the
global minimizer
is either the triangular bi-pyramid or else some
pyramid with square base.
The {\it triangular bi-pyramid\/} (TBP) is
the configuration of $5$ points with one
point at the north pole, one at the south
pole, and three arranged in an equilateral
triangle around the equator. 

[{\bf HS\/}] has a rigorous computer-assisted
proof that the TBP is the unique minimizer
for $R_{-1}$ (Polya's problem) and my paper [{\bf S1\/}] has a
rigorous computer-assisted proof that the TBP
is the unique minimizer for $R_1$ (Thomson's problem)
and $R_2$.
The paper [{\bf DLT\/}] gives a traditional
proof that the TBP is the unique minimizer
for the logarithmic potential.
The TBP is \underline{not} the minimizer
for $R_s$ when $s> 15.040809...$. 

Define
\begin{equation}
G_k(r)=(4-r^2)^k, \hskip 30 pt
k=1,2,3,...
\end{equation}
In [{\bf T\/}], A. Tumanov gives a traditional
proof of the following result.
\begin{theorem}[Tumanov]
\label{tumanov2X}
Let $f=a_1G_1+a_2G_2$ with $a_1,a_2>0$.  The
TBP is the unique global minimizer with respect to $f$.
Moreover, a critical point of $f$ must be the TBP
or a pyramid with square base.
\end{theorem}

As an immediate corollary, the TBP is a minimizer
for $G_1$ and $G_2$.
Tumanov points out that these potentials do not have
an obvious geometric interpretation, but they are
amenable to a tradicional analysis.  He points
out his result might be a step towards proving the that TBP
minimizes a range of power law potentials.
He makes the following observation:
If the TBP is the unique minimizer for $G_3$ and $G_5$,
then the TBP is the unique minimizer
for $R_s$ provided that $s \in (-2,2]-\{0\}$.

\subsection{Results}

Tumanov did not give a proof of his
observation in [{\bf T\/}], 
but we will prove a more extensive result.

\begin{lemma}
\label{Tumanov}
Suppose the TBP is the unique minimizer
for $G_3,G_4,G_5,G_6$ and
also for $G^{\#}_{10}=G_{10}+28G_5+102G_2.$
Then the TBP is the unique minimizer for
$R_s$ for any $s \in (-2,13]-\{0\}$.
\end{lemma}

\S \ref{general} has a clear description of
the idea behind the result.  This kind of
technique, in some form, is used in
many papers on energy minimization.  See
[{\bf BDHSS\/}] for very general ideas like this.

The main purpose of this paper is to prove
the following result.

\begin{theorem}[Main]
\label{main}
The TBP is the unique minimizer for 
$G_3,G_4,G_5,G_6$ and $G^{\#}_{10}$.
\end{theorem}

Combining Theorem \ref{main} with Lemma
\ref{Tumanov}, we get
\begin{theorem}
\label{aux}
The TBP is the unique
minimizer for $R_s$ if
$s \in (-2,13]-\{0\}$. 
\end{theorem}

The lower bound of $-2$ in Corollary \ref{aux} is
sharp, because the TBP is not a minimizer for $R_s$
when $s<-2$. We separate this case out for special
mention.

\begin{theorem}
\label{aux2}  Let $p \in (0,2)$ be arbitrary.
Then the TBP is the unique maximizer, amongst
$5$-point configurations 
on the unit sphere, of the sums of the
$p$th powers of the distances between the points.
\end{theorem}
In particular, our methods give another solution
of Polya's problem. This is the case $p=1$, which is
solved in [{\bf HS\/}] by a different kind of
computer-assisted proof.

\subsection{Outline of the Paper}

In \S 2 I will discuss the moduli space of
normalized configurations of $5$ points
on the sphere.  The key idea is to use stereographic
projection to move the points into $\R^2 \cup \infty$.
This gives the moduli space a natural Euclidean
structure, making a
divide-and-conquer algorithm easier to manage,
even though the expressions for the energy potentials
are more complicated.
Our basic object is a {\it block\/}, a
certain kind of rectangular solid in
the moduli space.

In \S 3 I will prove some elementary
facts about spherical geometry, and
define the main quantities associated to
blocks.  These quantities form the basis
of the main technical result, the Energy Theorem.

In \S 4 I will prove the Energy Theorem,
Theorem \ref{ENERGY}, which gives an
efficient estimate on the minimum energy
of all configurations contained within a
block $B$ based on the minimum energy of the
configurations associated to the vertices of
$B$ and an auxiliary error term.

In \S 5 I will describe a
divide-and-conquer algorithm based
on the Energy Theorem.  Given a
potential function $F$ from the Main Theorem
and a small neighborhood $B_0$ of the TBP, the 
program does a recursive search
through the poset of dyadic blocks, eliminating
a block $B$ if $B \subset B_0$ or if
the Energy Theorem determines that
all configurations in $B$ have higher 
energy than the TBP.  (Sometimes we also
eliminate $B$ by symmetry considerations.)
If $B$ is not eliminated,
$B$ is subdivided and the pieces are
then tested. If the program runs to
completion, it constitutes a proof that some
$F$-minimizer lies in $B_0$.
See Lemmas \ref{energy1} and
\ref{energy2}.

In \S 6 I prove that the Hessian of the
$F$ energy functional is positive definite
throughout $B_0$ for each relevant choice of $F$.
This combines with Lemmas
\ref{energy1} and \ref{energy2}
to prove the Main Theorem.
Our local analysis is essentially
Taylor's Theorem with Remainder.

In \S 7-8 I prove Lemma \ref{Tumanov}.

In \S 9 I will explain some details
about the computer program.  All
the calculations which go into
Lemmas \ref{energy1} and \ref{energy2}
are done with interval arithmetic
to avoid roundoff error.
The calculations for Lemma \ref{Tumanov}
are all exact integer calculations.
The only place where 
floating point calculations occur without
any control is for the
extremely loose bound in Equation \ref{GREAT}.
For convenience, we rely on
Mathematica [{\bf W\/}] to estimate some
numbers which involve square roots
of integers.  Whereas Mathematica can produce
decimal expansions of these numbers up to 
thousands of digits, we just need the 
answer to be accurate up to a factor of $3$.

\subsection{Discussion}

Our proof of Theorem \ref{main} will not work
for $G_1$ because the TBP is not the unique
minimizer in this case.
However, this case is relatively trivial.
Our proof would work for $G_2$, but the setup in
\S 2 needs to be modified somewhat. We don't
pursue this.

One can check that some pyramid with square base
has lower $G_7$ energy than the TBP.
I checked by calculating the hessian of
the energy, as in \S 6, 
that the TBP is not even a local minimizer
for $G_8,...,G_{100}$.  This fact led me to
search for $G_{10}^{\#}$ and other combinations like it.

The analysis of $G_{10}^{\#}$ is tougher than the
analysis of $G_3,G_4,G_5,G_6$.  For the (dis)interested
reader, we mention that one can still prove 
Theorem \ref{aux} for all $s \in (-2,6]-\{0\}$ without
$G_{10}^{\#}$.  Indeed, an earlier version of the
paper had this more limited result and then I decided
to push harder on the method.

I am starting to think that I can get the definitive
result on the positive side -- i.e. all the way to
the presumed cutoff of 15.040809...  Here is a sketch.
The same techniques we
use in \S 8 establish the following result.
\begin{lemma}
\label{Tumanov2}
Let $X_0$ denote the TBP
and let $X$ be any other $5$-point configuration on the
unit sphere.  Suppose that the energy of $X$ exceeds
the energy of $X_0$ with respect to
$$G_2^*=3G_2-2G_1, \hskip 20 pt
G_5^*=G_5-5G_1, \hskip 30 pt
G_{10}^*=G_{10}+13 G_5 + 55 G_2.$$
Then $X$ is not the minimizer
for the $R_s$ energy when $s \in [13,15.05]$.
\end{lemma}
My program runs to completion
on $G_2^*$ and $G_5^*$.  Since the TBP is not
the minimizer for $G_{10}^*$, my program will
not run to completion. However, the program
eliminates most other configurations.
Precisely, my program (coupled with
the same local analysis as is done in
\S 6) seems to show that a minimizer
for $G_{10}^*$ either is the TBP or else
can be rotated so that its stereographic
projection has points $z_0,z_1,z_2,z_3,\infty$,
where
\begin{itemize}
\item $z_0 \in [13/16] \times \{0\}$.
\item $z_1 \in [-1/16,1/16] \times [-5/8,-15/16]$.
\item $z_2 \in [-1,-3/4] \times [-1/16,1/16]$.
\item $z_3 \in [-1/16,1/16] \times [5/8,15/16]$.
\end{itemize}
Here is a picture of these regions in the plane.
The $4$ dots represent the stereographic
projection of the nearest copy of the TBP.
The grey circle is the unit circle.

\begin{center}
\resizebox{!}{2.8in}{\includegraphics{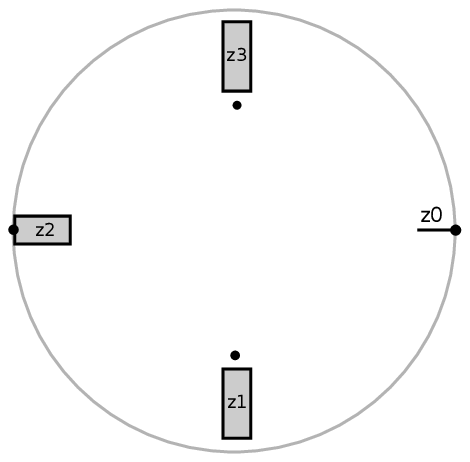}}
\newline
{\bf Figure 2:\/} Location of the minimizer for
power laws in $[13,15.05]$.
\end{center}

These configurations are quite close to pyramids
with square base, and for all of them it
appears (from a billion random trials)
that the symmetrization
operation
$$(z_0,z_1,z_2,z_3) \to (\rho_0,-i \rho_1,-\rho_0,i\rho_1),
\hskip 30 pt \rho_i=\bigg|\frac{z_i-z_{i+1}}{2}\bigg|$$
decreases the energy with respect to $G_k$ for
all $k=2,...,14$.  This result, combined with another
approximation lemma like Lemma \ref{Tumanov2}, seems
sufficient to show that for any $s \in [13,15.05]$
the minimizer with respect to $R_s$ is either
the TBP or a pyramid with square base.  That would
show the existence of a cutoff around $15.040809$ and allow
one to compute its value, say, to a million decimal places.
\newline

This paper has many points in common with [{\bf S1\/}] but
the exposition here is simpler and cleaner.   Now I have a 
better understanding of which ingredients are important
in the running time of the program, I was able to
give a more efficient treatment of the problem.  In
the end, the estimates and the routine are pretty simple.
I think that a competent programmer could reproduce the
results.

As in [{\bf S1\/}], the crucial feature
of this paper is a good
error estimate, the Energy Theorem, which allows one to control
the energies of all configurations within a
block of the configuration space just by considering
finitely many of those configurations.  A good
estimate means the difference between
a feasible calculation and one that would outlast
the universe.  I view the Energy Theorem as the
main mathematical contribution of the paper.

All the estimates that
go into the Energy Theorem are rational functions
of the inputs.  So, one could in theory
give a proof that just uses integer arithmetic.
An early version of my program ran integer
arithmetic calculations, but they seemed very
slow.  I hope to eventually optimize the
code for integer calculations and try it out
seriously.  It would be nice to have an
exact integer proof of the Main Theorem.

\subsection{Companion Computer Programs}

The computer code involved in this paper is
publicly available by download from my website.
See \S \ref{download} for details.
There is a small amount of Mathematica code,
and then there is are two
large Java programs.  The main java program
runs the divide-and-conquer program which
establishes Lemma \ref{energy1} and \ref{energy2}.
The other java program does the calculations
for Lemma \ref{Tumanov}.  The Mathematica code
is used in the local analysis of the Hessian
in \S 6.  The Java programs come with
graphical user interfaces, and each one has built
in documentation and debuggers. The interfaces
let the user watch the programs in action,
and check in many ways that they are operating
correctly.

\subsection{Acknowledgements}

I would like to thank Henri Cohn, Abhinav Kumar,
Curtis McMullen, Jill Pipher, and Sergei Tabachnikov
for discussions related to this paper. I would like
to thank Alexander Tumanov for his great observation.
Also, I would like to thank I.C.E.R.M. for facilitating
this project. My interest in this problem
was rekindled during discussions about point configuration
problems around the time of the I.C.E.R.M. Science
Advisory Board meeting in Nov. 2015.

\newpage

%% file: 2config.tex
\section{The Configuration Space}

\subsection{Normalized Configurations}
\label{pn}

Recall that
\begin{equation}
G_k(r)=(4-r^2)^k, \hskip 30 pt
G^{\#}_{10}=G_{10}+28 G_5 + 102 G_2.
\end{equation}
The TBP has $6$ bonds of length $\sqrt 2$,
and $3$ bonds of length $\sqrt 3$, and one
bond of length $2$.  Hence, the $G_k$
energy of the TBP is $3+6\times 2^k$.  In particular,
\begin{equation}
\label{energy0}
({\cal E\/}_3,{\cal E\/}_4,{\cal E\/}_5,{\cal E\/}_6,
{\cal E\/}^{\#}_{10})=
(51,99,195,387,14361).
\end{equation}

Let $\widehat p_0,...,\widehat p_4$ be a configuration of
$5$ points on $S^2$.  We call this
configuration {\it normalized\/} if
\begin{equation}
\widehat p_4=(0,0,1), \hskip 30 pt
\|\widehat p_4-\widehat p_0\| \leq
\|\widehat p_4-\widehat p_i\|, \hskip 20 pt
\forall i \in \{1,2,3\}.
\end{equation}
The choice of the bizarre constant in the next lemma
will become clearer in the next section.

\begin{lemma}
\label{config1}
For a normalized minimizer, we have
$\|\widehat p_4-\widehat p_0\|>1/2$ and
$\|\widehat p_4-\widehat p_i\|>4/\sqrt{13}$
for $i=1,2,3$.
\end{lemma}

\startproof
The values of $G_3,G_4,G_5,G_6,G^{\#}_{10}$
evaluated at $1/2$ respectively exceed
the values in Equation \ref{energy0}, so a
normalized minimizer cannot have a point within
$1/2$ of $\widehat p_4$.

If $\|\widehat p_4-\widehat p_i\| \leq 4/\sqrt{13}$ for some
$i=1,2,3$ then $\|\widehat p_4-\widehat p_0\| \leq 4/\sqrt{13}$ as well.
But then
the sum of the energies coming from 
bonds between $\widehat p_4$
and other points is at least
$2G_k(4/\sqrt 13)=2(36/13)^k$. This
exceeds the value given in
Equation \ref{energy0} except in
the case of $G_3$, where all we get
is $2(36/13)^3>42$.

We need a special argument for $G_3$.
Notice that the same argument as above
shows that the distance between
any two points for the $4$-point minimizer is
at least $1$, because otherwise the energy of the
single bond, $3^k$, would exceed the energy $6 \times (4/3)^k$
of the regular tetrahedron.
Now observe that the function $G_3$ is convex
decreasing on the interval $[1,2]$.  But the
regular tetrahedron is the mimimizer
for any convex decreasing potential -- see
[{\bf CK\/}] for a proof.  Hence, 
the regular tetrahedron
is the global minimizer for $G_3$ amongst all
$4$ point configurations.   Hence
the sum of energies in the bonds not involving
$\widehat p_4$ is at least $6 \times (4/3)^3>14$.
Since $14+42>51$,
our configuration could not be an minimizer.
\endproof

\subsection{Stereographic Projection}

As in [{\bf S1\/}] we work mainly with
$\R^2 \cup \infty$ rather than on $S^2$.    
Our reason for this is that a configuration space
based on points in $\R^2$ has a
natural flat structure, and lends itself well
to divide-and-conquer algorithms.

We map $S^2$ to $\R^2 \cup \infty$ using
{\it stereographic projection\/}:
\begin{equation}
\label{stereo}
\Sigma(x,y,z)=\bigg(\frac{x}{1-z}, \frac{y}{1-z}\bigg).
\end{equation}
$\Sigma$ is a conformal diffeomorphism which maps circles
in $S^2$ to lines and circles in $\R^2 \cup \infty$.
The inverse is:
\begin{equation}
\label{inversestereo}
\Sigma^{-1}(x,y)=\bigg(\frac{2x}{1+x^2+y^2},\frac{2y}{1+x^2+y^2},
1-\frac{2}{1+x^2+y^2}\bigg).
\end{equation}

We will use the convention that any object $S$ in $\R^2 \cup \infty$
corresponds to 
\begin{equation}
\widehat S=\Sigma^{-1}(S)
\end{equation}
in $S^2$.   Given our configuration
$\widehat p_0,...,\widehat p_4$ as in
the previous section, we have points
$p_0,...,p_3 \in \R^2$ and $p_4=\infty$.
Thus, the points $p_0,...,p_3$ determine
the configation.

\begin{lemma}
If $p_0,p_1,p_2,p_3$ correspond to a normalized
minimizer for $G_k$ with $k \in \{3,4,5,6\}$, then
\label{farout}
$|p_0|<4$ and
$|p_i|<3/2$ for $i=1,2,3$.
\end{lemma}

\startproof
We compute that
\begin{equation}
\|\Sigma^{-1}(4,0)-(0,0,1)\|<1/2, \hskip 30 pt
\|\Sigma^{-1}(3/2,0)-(0,0,1)\|=4/\sqrt{13}.
\end{equation}
By symmetry and monotonicity,
the first equation holds for any $p$ with $|p| \geq 4$
and the second equation holds for any $p$ with
$|p| \geq 3/2$.  Lemma \ref{config1}
finishes the proof.
\endproof

By symmetry we can rotate the picture so that
$p_0$ lies in the $x$-axis, and has coordinate
$(p_{01},0)$, with $p_{01} \in (0,4)$.
The remaining points must lie in the
disk of radius $3/2$ about $0$.  However, we
find it convenient to confine $p_1,p_2,p_3$ only to
the square $[-2,2]^2$ at first.
We set our moduli space to
be the product
\begin{equation}
\square= [0,4] \times \Big([-2,2]^2\Big)^3=
[0,4] \times [-2,2]^6.
\end{equation}
Conveniently,
$\square$ is a cube of sidelength $4$
in $\R^7$.  We just have to search 
through this big cube and eliminate everything
which is not the TBP!

\subsection{The TBP Configurations}
\label{TBP}

The TBP has two kinds of points, the two at the
poles and the three at the equator.
When $\infty$ is a polar point, the points
$p_0,p_1,p_2,p_3$ are, after suitable permutation,
\begin{equation}
(1,0), \hskip 20 pt
(-1/2,-\sqrt 3/2), \hskip 20 pt
(0,0), \hskip 20 pt
(-1/2,+\sqrt 3/2).
\end{equation}
We call this the {\it polar configuration\/}.
When $\infty$ is an equatorial point, the points
$p_0,p_1,p_2,p_3$ are, after suitable permutation,
\begin{equation}
(1,0), \hskip 20 pt (0,-1/\sqrt 3),
\hskip 20 pt
(-1,0), \hskip 20 pt
1/\sqrt 3).
\end{equation}
We call this the {\it equatorial configuration\/}.
We prefer the polar configuration, and we will use
a trick to avoid ever having to search very near the
equatorial configuration.

We can visualize the two configurations together 
in relation to the regular $6$-sided star. The black
points are part of the polar configuration and
the white points are part of the equatorial configuration.
The grey point belongs to both configurations.
The points represented by little squares are polar 
and the points represented by little disks are equatorial.
The beautiful pattern made by these two configurations
is not part of our proof, but it is nice to contemplate.

\begin{center}
\resizebox{!}{2.8in}{\includegraphics{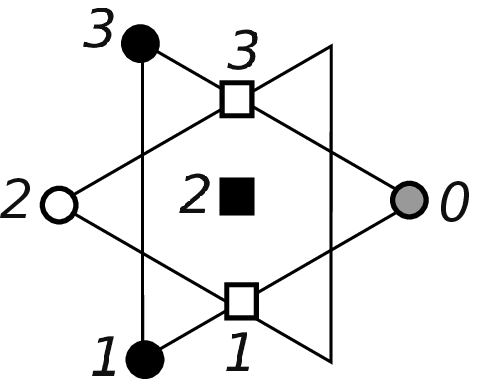}}
\newline
{\bf Figure 2.1:\/} Polar and equatorial versions of the
TBP.
\end{center}

\subsection{Dyadic Blocks}
\label{DYAD}

In $1$ dimension, the {\it dyadic subdivision\/} of a
line segment is the union of the two segments obtained
by cutting it in half.  In $2$ dimensions, the
{\it dyadic subdivision\/} of a square is the
union of the $4$ quarters that result in cutting
the the square in half along both directions.
See Figure 2.2.

\begin{center}
\resizebox{!}{.7in}{\includegraphics{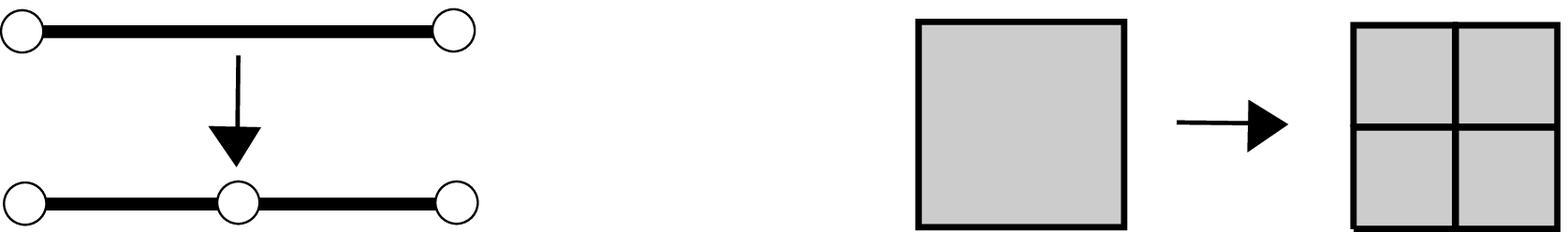}}
\newline
{\bf Figure 2.2:\/} Dyadic subdivision
\end{center}

We say that a {\it dyadic segment\/} is any segment
obtained from $[0,4]$ by applying dyadic subdivision
recursively.  We say that a {\it dyadic square\/} 
is any square obtained from $[-2,2]^2$ by applying
dyadic subdivision recursively.
We count $[0,4]$ as a dyadic segment and
$[-2,2]^2$ as a dyadic square.
\newline
\newline
{\bf Hat and Hull Notation:\/}
Since we are going to be switching back and forth between
the picture on the sphere and the picture in $\R^2$, we
want to be clear about when we are talking about
solid bodies, so to speak, and when we are talking about
finite sets of points.
We let $\langle X \rangle$ denote the convex hull
of any Euclidean subset. Thus, we think of a dyadic
square $Q$ as the set of its $4$ vertices and
we think of $\langle Q \rangle$ as the
solid square having $Q$ as its vertex set.
Combining this with our notation for stereographic
projection, we get the following notation:
\begin{itemize}
\item $\widehat Q$ is a set of $4$ co-circular points
on $S^2$.
\item $\langle \widehat Q \rangle$ is a convex quadrilateral
whose vertices are $\widehat Q$.
\item $\widehat{\langle Q \rangle}$ is a ``spherical patch''
on $S^2$, bounded by $4$ circular arcs.
\end{itemize}
We will use this notation throughout the paper.
\newline
\newline
\noindent
{\bf Good Squares:\/}
A dyadic square is {\it good\/} if it
is contained in $[-3/2,3/2]^2$ and has
side length at most $1/2$.
Note that a good dyadic square cannot cross the 
coordinate axes.  The only dyadic square which
crosses the coordinate axes is $[-2.2]^2$, and
this square is not good.
Our computer program will only do spherical
geometry calculations on good squares.
\newline
\newline
{\bf Dyadic Blocks:\/}
We define a {\it dyadic block\/} to be a
$4$-tuple $(Q_0,Q_1,Q_2,Q_3)$, where
$Q_0$ is a dyadic segment and $Q_i$ is a
dyadic square for $j=1,2,3$.  We say that
a block is {\it good\/} if each of its $3$ component
squares is good.
By Lemma \ref{farout}, any energy minimizer for
$G_k$ is contained in a good block. Our
algorithm in \S \ref{ALG} 
quickly chops up the blocks in
$\square$ so that only good ones are considered.

The product
\begin{equation}
\langle B \rangle=\langle Q_0 \rangle \times
\langle Q_1 \rangle \times
\langle Q_2 \rangle \times
\langle Q_3 \rangle
\end{equation}
is a rectangular solid in the
configuration space $\square$.
On the other hand, the product
\begin{equation}
B= Q_0  \times
 Q_1  \times
 Q_2  \times
 Q_3 
\end{equation}
is the collection of $128$ vertices of
$\langle B \rangle$. 
We call these the
{\it vertex configurations\/} of the block.
\newline
\newline
{\bf Definition:\/}
We say that a configuration
$p_0,p_1,p_2,p_3$ is {\it in\/} the block $B$
if $p_i \in \langle Q_i \rangle$ for $i=0,1,2,3$.
In other words, the point in $\square$ representing
our configuration is contained in $\langle B \rangle$.
Sometimes we will say that this configuration is
{\it associated to\/} the block.
\newline
\newline
{\bf Sudvidision of Blocks:\/}
There are $4$ obvious subdivision operations we can
perform on a block. 
\begin{itemize}
\item The operation $S_0$ divides
$B$ into the two blocks
$(Q_{00},Q_1,Q_2,Q_3)$ and
$(Q_{01},Q_1,Q_2,Q_3)$.
Here $(Q_{00},Q_{01})$ is the dyadic subdivision of $Q_0$.
\item the operation $S_1$ divides
$B$ into the $4$ blocks
$(Q_0,Q_{1ab},Q_2,Q_3)$, where
$(Q_{100},Q_{101},Q_{110},Q_{111})$ is
the dyadic subdivision of $Q_1$.
\end{itemize}
The operations $S_2$ and $S_3$ are similar
to $S_1$.

The set of dyadic blocks has a natural poset structure
to it, and basically our algorithm does a depth-first-search
through this poset, eliminating solid blocks either
due to symmetry considerations or due to energy 
considerations.  The short chapter \S \ref{ALG}
describes the algorithm precisely.

\subsection{Relevant Configurations}
\label{normalized}
\label{checkSP}

We write $p_i=(p_{i1},p_{i2})$.
Call $(p_0,p_1,p_2,p_3) \in \square$
{\it relevant\/} if all $5$ of these inequalities hold:
\begin{equation}
p_{01} \geq 1, \hskip 15 pt
p_{12} \leq 0, \hskip 15 pt
p_{22} \geq 0, \hskip 15 pt
p_{22} \leq p_{32}, \hskip 15 pt
p_{21} \geq -1/2.
\end{equation}
and otherwise {\it irrelevant\/}.
The polar TBP is relevant but
the equatorial TBP is irrelevant.

\begin{lemma}
Relative to any monotone decreasing energy potential,
some minimal configuration is relevant.
\end{lemma}

\startproof
We normalize so that $\widehat p_0$ and
$\widehat p_4$ are the closest pair of points.
If Inequality 1 fails then $\widehat p_0$ and
$\widehat p_4$ are more than $\sqrt 2$ apart,
and hence all the points are.
But it is impossible to place $5$
points on $S^2$ such that every two of them are more than
$\sqrt 2$ apart.  See [{\bf S1\/}] for a proof.

If $p_1,p_2,p_3$ lie all inside (or all outside) 
the upper half plane $H$ then
all $5$ corresponding points in $S^2$ lie in the
hemisphere $\widehat H$ 
and $\widehat p_1,\widehat p_2,\widehat p_3$
lie in the interior of $\widehat H$.
But then we can get a configuration of lower
energy by replacing $\widehat p_1$ with its reflection in
$\partial \widehat H$.  
So, reflecting in the $x$-axis and relabeling if
necessary we can assume that $p_1$ does not lie in
the interior of $H$ and $p_2,p_3$ lie in $H$.
This gives Inequalities 2,3,4.

If Inequality 5 fails, we apply
the inversion $I_C$ in the circle $C$ which is
centered at $p_0$ and has the property that
$\Sigma^{-1}(C)$ is a great circle. $I_C$
corresponds to the isometric reflection in 
$\Sigma^{-1}(C)$ which swaps $\widehat p_0$
and $\widehat p_4$.  Let $q_j=I_C(p_j)$ for $j=1,2,3$.
Let $\Delta$ be the disk bounded by $C$.
Note that $I_C$ fixes $C$ and swaps $\Delta$ with
the complment of $\Delta$.
Let $X$ denote the $x$-axis.
If $p_0=(1,0)$, then $C \cap X=\{1 \pm \sqrt 2\}$.
This is the extreme case. When $p_{01}>1$ both 
points of $C \cap X$ lies to the right of
$(1-\sqrt 2,0)$, which in
turn lies to the right of $(-1/2,0)$.  But then
$p_2 \not \in \Delta$.  Hence
$q_2 \in \Delta$.  This gives $q_{21}>-1/2$.

Since $I_C$ preserves all
the rays through $p_0$, we still have
$q_{12} \leq 0$ and $q_{22},q_{32} \geq 0$.
If $q_3 \in \Delta$ then $q_{31}>-1/2$.  If needed,
we switch $q_2$
and $q_3$ to guarantee that $q_{22} \leq q_{32}$.
If $q_3 \not \in \Delta$ then $p_{31}>1-\sqrt 2>p_{21}$ and hence
the ray $\rho_3$ joining $p_0$ to $p_3$ lies above the
ray $\rho_2$ joining $p_0$ to $p_2$. But then
$q_3 \in \rho_3-\Delta$ and $q_2 \in \rho_2 \cap \Delta$.
This forces $q_{22}<q_{32}$.  
So in all cases, we can retain Inequalities 2,3,4.
Since $q_0=p_0$ we retain Inequality 1.
\endproof

Call a block {\it irrelevant\/} if every
configuration in the interior of the
block is irrelevant.  Call a block
{\it relevant\/} if it is not irrelevant.
Now we give a criterion for a block to be
irrelevant.

Given a box $Q_j$, let $\overline Q_{jk}$ and
$\underline Q_{jk}$
denote the maximum and
minimum $k$th coordinate of a point in $Q_j$.
The good block $B=(Q_0,Q_1,Q_2,Q_3)$ is irrelevant
provided that at least one of the following
holds.
\begin{equation}
\label{irrelevantX}
\overline Q_{01} \leq 1, \hskip 15 pt
\underline Q_{12} \geq 0, \hskip 15 pt
\overline Q_{22} \leq 0, \hskip 15 pt
\underline Q_{22}\geq \overline Q_{32}, \hskip 15 pt
\overline Q_{21} \leq -1/2.
\end{equation}
Notice that these $5$ inequalities parallel
the ones given above for individual configurations.

Every relevant configuration in
the boundary of an irrelevant block is
also in the boundary of a relevant block.
So, to prove Theorem \ref{main}, we can
ignore the irrelevant blocks.

\subsection{A Technical Result about Dyadic Boxes}

The following result will be useful for
the estimates in \S \ref{sep}.  We state it
in more generality than we need, to show
what hypotheses are required, but we note
that good dyadic squares
satisfy the hypotheses.  We only
care about the result for
good dyadic squares and for
good dyadic segments. In the
case of good dyadic segments, the
lemma is obviously true.

\begin{lemma}
\label{mono}
Let $\widehat Q$ be a rectangle whose
sides are parallel to the coordinate
axes and do not cross the coordinate
axes.  Then the points of
$\widehat{\langle Q \rangle}$ closest
to $(0,0,1)$ and farthest from $(0,0,1)$ are
both vertices.
\end{lemma}

\startproof
Put the metric on $\R^2 \cup \infty$ which makes
stereographic projection an isometry.  By symmetry,
the metric balls about $\infty$ are the complements
of disks centered at $0$.  The smallest disk centered
at $0$ and containing $\langle Q \rangle$ must have a
vertex on its boundary.  Likewise, the largest disk
centered at $0$ and disjoint from the interior
of $\langle Q \rangle$ must have a vertex in its
boundary.  This latter result uses the fact that
$\langle Q \rangle$ does not cross the
coordinate axes.  These statements are equivalent
to the statement of the lemma.
\endproof

\subsection{Very Near the TBP}
\label{checkTBP}

Our calculations will depend on a
a pair $(S,\epsilon_0)$, both powers of two.
For $G_3,G_4,G_5,G_6$ we use
$S=2^{25}$ and $\epsilon_0=2^{-15}$.
For $G^{\#}_{10}$ we use
$S=2^{30}$ and $\epsilon_0=2^{-18}$.
Our discussion below would work for 
other similar choices.

Define the {\it in-radius\/} of a cube
to be half its side length.  Let
$P_0$ denote the configuration representing the
totally normalized polar TBP and let
$B_0$ denote the cube centered at $P_0$
ahd having in-radius $\epsilon_0$.
Note that $B_0$ is not a dyadic
block.  This does not bother us. 
Here we give a sufficient condition for
$B \subset B_0$, where $B$ is some dyadic block.

Let $a=\sqrt 3/2$.  For each choice of $S$ we
compute a value $a^*$ such that
$Sa^* \in \Z$ and $|a-a^*|<1/S$.  There are
two such choices, namely
\begin{equation}
\frac{{\rm floor\/}(Sa)}{S}, \hskip 30 pt
\frac{{\rm floor\/}(Sa+1)}{S}.
\end{equation}  
In practice, our program sets $a^*$ to be the
first of these two numbers, but we want to
state things in a symmetric
way that works for either choice.

We define
$B'_0=Q'_0 \times Q'_1 \times Q'_2 \times Q'_3$ where
$$Q_0'=[1-\epsilon_0,1+\epsilon_0],$$
$$Q_1'=[-1/2-\epsilon_0,-1/2+\epsilon_0] 
\times [-S^{-1}-a^*-\epsilon_0,S^{-1}-a^*+\epsilon_0]$$
$$Q_2'=[-\epsilon_0,\epsilon_0]^2$$
$$Q_3'=[1/2-\epsilon_0,1/2+\epsilon_0] \times 
\times [-S^{-1}-a^*-\epsilon_0,S^{-1}-a^*+\epsilon_0]$$
By construction, we have
$B'_0 \subset B_0$.
Given a block $B=Q_0 \times Q_1 \times Q_2 \times Q_3$,
the condition
$Q_i \subset Q_i'$ for all $i$ implies
that $B \subset B_0$.   
We will see in \S \ref{computer} that this is an
exact integer calculation for us.

\newpage

%% file: 3sphere.tex
\section{Spherical Geometry Estimates}

\subsection{Overview}

In this chapter we define the basic quantities
that go into the Energy Theorem, Theorem
\ref{ENERGY}.
We will persistently use the hat and hull
notation defined in \S \ref{DYAD}.  Thus:
\begin{itemize}
\item When $Q$ is a dyadic square,
$\langle \widehat Q \rangle$ is a
convex quadrilateral in space whose vertices
are $4$ co-circular points on $S^2$, and
$\widehat{\langle Q \rangle}$ is a
subset of $S^2$ bounded by $4$ circular arcs.
\item When $Q$ is a dyadic segment, 
$\langle \widehat Q \rangle$ is a
segment whose endpoints are on $S^2$, and
$\widehat{\langle Q \rangle}$ is an arc of a great circle.
\end{itemize}

Here is a summary of the quantities we will
define in this chapter.  The first three
quantities are not rational functions of the
inputs, but our estimates only use the
squares of these quantities.
\newline
\newline
{\bf Hull Diameter:\/} $d(Q)$ will be the diameter
$\langle \widehat Q \rangle$.  
\newline
\newline
{\bf Edge Length:\/} $d_1(Q)$ will be the length
of the longest edge of 
$\langle \widehat Q \rangle$. 
\newline
\newline
{\bf Circular Measure:\/} Let $D_Q \subset \R^2$
denote the disk containing 
$Q$ in its boundary.  $d_2(Q)$ is the
diameter of $\widehat D_Q$. 
\newline
\newline
{\bf Hull Separation Constant:\/}
$\delta(Q)$ will be a constant such that every
point in
$\widehat{\langle Q \rangle}$ is within
$\delta(Q)$ of a point of
$\langle \widehat Q \rangle$.
This quantity is a rational function
of the coordinates of $Q$.
\newline
\newline
{\bf Dot Product Bounds:\/}
We will introduce a (finitely computable, rational)
 quantity
$(Q \cdot Q')_{\rm max\/}$ which
has the property that
$$V \cdot V' \leq
(Q \cdot Q')_{\rm max\/}$$ for all 
$V \in \widehat{\langle Q \rangle} \cup \langle \widehat Q \rangle$
and
$V' \in \widehat{\langle Q' \rangle} \cup \langle \widehat Q' \rangle$.

\subsection{Some Results about Circles}

Here we prove a few geometric facts about
circles and stereographic projection.

\begin{lemma}
\label{diam}
Let $Q$ be a good dyadic square or a dyadic
segment.  The circular arcs bounding
$\widehat{\langle Q \rangle}$ lie in
circles having diameter at least $1$.
\end{lemma}

\startproof
Let $\Sigma$ denote stereographic projection.
In the dyadic segment case, $\widehat{\langle Q \rangle}$
lies in a great circle.
In the good dyadic square case, each edge of $\langle Q \rangle$
lies on a line $L$ which contains a point $p$ at most $3/2$ from the origin.
But $\Sigma^{-1}(p)$ is at least $1$ unit from $(0,0,1)$.
Hence $\Sigma^{-1}(L)$, which limits on $(0,0,1)$ and
contains $\Sigma^{-1}(p)$,
has diameter at least $1$.  The set
$\Sigma^{-1}(L \cup \infty)$ is precisely
the circle extending the relevant edge of 
$\widehat{\langle Q \rangle}$
\endproof

Let $D \subset \R^2$ be a disk of radius $r \leq R$
centered at a point which is $R$ units
from the origin. 
Let $\widehat D$ denote the corresponding disk
on $S^2$.  We consider $\widehat D$ as a subset
of $\R^3$ and compute its diameter in with
respect to the Euclidean metric on $\R^3$.

\begin{lemma}
\begin{equation}
{\rm diam\/}^2(\widehat D)=
\frac{16r^2}{1+2r^2+2R^2+(R^2-r^2)^2}.
\end{equation}
\end{lemma}

\startproof
By symmetry it suffices to consider the
case when the center of $D$ is $(R,0)$.  The
diameter is then achieved by the two
points $V=\Sigma^{-1}(R-r,0)$ and
$W=\Sigma^{-1}(R+r,0)$. 
 The formula comes
from computing $\|V-W\|^2$ and simplifying.
\endproof

We introduce the functions
\begin{equation}
\chi(D,d)=\frac{d^2}{4D}+\frac{d^4}{2D^3}
\hskip 30 pt
\chi^*(D,d)=\frac{1}{2}(D-\sqrt{D^2-d^2}).
\end{equation}
The second of these is a function closely
related to the geometry of circles.
This is the function we would use if we had
an ideal computing machine.  However, since
we want our estimates to all be rational
functions of the inputs, we will use the
first function.  We first prove an approximation
lemma and then we get to the main point.

\begin{lemma}
If $0 \leq d \leq D$ then $\chi^*(D,d) \leq \chi(D,d)$.
\end{lemma}

\startproof
If we replace $(d,D)$ by $(rd,rD)$ then both sides
scale up by $r$.  Thanks to this homogeneity, it suffices
to prove the result when $D=1$.  We have
$\chi(1,1)=3/4>1/2=\chi^*(1,1)$.  So, if suffices
to prove that the equation
$$\chi(1,d)-\chi^*(1,d)=\frac{d^2}{4}+\frac{d^4}{2} - 
\frac{1}{2}(1-\sqrt{1-d^2})$$
has no real solutions in $[0,1]$ besides $d=0$.
Consider a solution to this equation.
Rearranging the equation, we get $A=B$ where
$$A=-\frac{1}{2}\sqrt{1-d^2},
\hskip 30 pt B=\frac{d^2}{4}+\frac{d^4}{2} - 1/2.$$
An exercise in calculus shows that the
only roots of $A^2-B^2$ in $[0,1]$ are $0$ and
$\sqrt{1/2(\sqrt 8-1)}>.95$.  On the other hand
$A<0$ and $B>0$ on $[.95,1]$.
\endproof

Now we get to the key result.  This result holds in
any dimension, but we will apply it once to the
$2$-sphere, and once to circles contained in suitable
planes in $\R^3$.

\begin{lemma}
\label{circle}
Let $\Gamma$ be a round sphere of diameter $D$, contained
in some Euclidean space. Let $B$ be the
ball bounded by $\Gamma$.
Let $\Pi$ be a hyperplane which intersects
$B$ but does
not contain the center of $B$.
Let $\gamma=\Pi \cap B$ and let
$\gamma^*$ be the smaller of the
two spherical caps on $\Gamma$
bounded by $\Pi \cap \Gamma$.
Let $p^* \in \gamma^*$ be a point.
Let $p \in \gamma$ be the point so
the line $\overline{pp^*}$ contains
the center of $B$. Then
$\|p-p^*\| \leq  \chi(D,d)$.
\end{lemma}

\startproof
The given distance is maximized when
$p^*$ is the center of $\gamma^*$ and
$p$ is the center of $\gamma$. 
In this case it suffices by symmetry
to consider the situation in $\R^2$,
where $\overline{pp^*}$ is the perpendicular
bisector of $\gamma$.
Setting $x=\|p-p^*\|$, we have
\begin{equation}
\label{cross}
x(D-x)=(d/2)^2.
\end{equation}
This equation comes from a well-known theorem
from high school geometry concerning the lengths
of crossing chords inside a circle.  When we solve
Equation \ref{cross} for $x$, we see that
$x=\chi^*(D,d)$.  The previous lemma finishes the proof.
\endproof

\subsection{The Hull Approximation Lemma}
\label{hullsep}

\noindent
{\bf Circular Measure:\/}
When $Q$ is a dyadic square or segment, we
define
\begin{equation}
d_2(Q)={\rm diam\/}(\widehat D_Q),
\end{equation}
Where $D_Q \subset \R^2$ is such that
$Q \subset D_Q$.
So $d_2(Q)$ is the diameter of the small spherical cap
which contains $\widehat Q$ in its boundary.
Note that $\widehat{\langle Q \rangle} \subset \widehat D_Q$
by construction.  We call $d_2(Q)$ the {\it circular measure\/} of $Q$.
\newline
\newline
{\bf Hull Separation Constant:\/}
Recall that $d_1(Q)$ is the maximum
side length of $\langle Q \rangle$.  When
When $Q$ is a dyadic segment, we define
$\delta(Q)=\chi(2,d_2)$. 
When $Q$ is a good dyadic
square,
We define
\begin{equation}
\delta(Q)=\max\Big(\chi(1,d_1), \chi(2,d_2)\Big).
\end{equation} 
This definition makes sense, because
$d_1(Q) \leq 1$ and
$d_2(Q) \leq \sqrt 2<2$.  The point
here is that $\Sigma^{-1}$ is $2$-Lipschitz
and $Q$ has side length at most $1/2$.
We call $\delta(Q)$ the 
{\it Hull approximation constant\/}
of $Q$.

\begin{lemma}[Hull Approximation]
Let $Q$ be a dyadic segment or a
good dyadic square.
Every point of the spherical patch
$\widehat{\langle Q \rangle}$ is within
$\delta(Q)$ of a point of the
convex quadrilateral
$\langle \widehat Q \rangle$.
\end{lemma}

\startproof
Suppose first that $Q$ is a dyadic segment.
$\widehat{\langle Q \rangle}$ is
the short arc of a great circle
sharing endpoints with
$\langle \widehat Q \rangle$, a chord
of length $d_2$.  
By Lemma \ref{circle} each
point of on the circular arc is
within $\chi(2,d_2)$ of a
point on the chord.

Now suppose that $Q$ is a good dyadic square.
Let $O$ be the origin in $\R^3$.
Let $H \subset S^2$ denote
the set of points such that the segment
$Op^* \in H$ intesects
$\langle \widehat Q \rangle$ in a point $p$.  
Here $H$ is the itersection with the cone
over $\langle \widehat Q \rangle$ with $S^2$.
\newline
\newline
{\bf Case 1:\/}
Let $p^* \in \widehat{\langle Q \rangle} \cap H$.
Let $p \in \langle Q \rangle $ be such that
the segment $Op^*$ contains $p$. Let $B$
be the unit ball. Let
$\Pi$ be the plane containing $\widehat Q$.
Note that $\Pi \cap S^2$ bounds the spherical 
$\widehat D_Q$ which contains
$\widehat Q$ in its boundary. Let
$$
\Gamma=S^2, \hskip 20 pt
\gamma^*=\widehat D_Q, \hskip 20 pt
\gamma=\Pi \cap B.$$
The diameter of $\Gamma$ is $D=2$.  Lemma
\ref{circle} now tells us that
$\|p-p^*\| \leq \chi(2,d_2)$.
\newline
\newline
{\bf Case 2:\/}
Let $p^*\in \widehat{\langle Q\rangle}-H$.
The sets $\widehat{\langle Q\rangle}$ and
$H$ are both bounded by $4$ circular
arcs which have the same vertices. $H$
is bounded by arcs 
of great circles and $\widehat{\langle Q \rangle}$ is
bounded by arcs 
of circles having diameter at least $1$.
The point $p^*$ lies between an edge-arc $\alpha_1$ of
$H$ and an edge-arc $\alpha_2$ of $\widehat{\langle Q \rangle}$ 
which share both
endpoints.  Let $\gamma$ be the line segment
joining these endpoints.  The diameter of
$\gamma$ is at most $d_1$.

Call an arc of a circle {\it nice\/} if it is
contained in a semicircle, and if the circle
containing it has diameter at least $1$.
The arcs $\alpha_1$ and $\alpha_2$ are both nice.
We can foliate the region between $\alpha_1$ and
$\alpha_2$ by arcs of circles. These circles are
all contained in the intersection of $S^2$ with
planes which contain $\gamma$.  Call this
foliation $\cal F$.   We get this foliation by
rotating the planes around their common
axis, which is the line through $\gamma$.

Say that an ${\cal F\/}$-{\it circle\/} is
a circle containing an arc of $\cal F$.
Let $e$ be the edge of $\langle Q \rangle$
corresponding to $\gamma$.  
Call $(e,Q)$ {\it normal\/} if
$\gamma$ is never the diameter of an
$\cal F$-circle.  If $(e,Q)$  is normal, 
then the diameters of
the $\cal F$-circles interpolate
monotonically between the diameter of
$\alpha_1$ and the diameter of $\alpha_2$.
Hence, all $\cal F$-circles have diameter
at least $1$.  At the same time, if
$(e,Q)$   is normal, then all arcs of
$\cal F$ are contained in semicircles,
by continuity.  In short, if $(e,Q)$ is normal,
then all arcs of $\cal F$ are nice.
Assuming that $(e,Q)$ is normal,
let $\gamma^*$ be the arc in $\cal F$
which contains $p^*$. Let $p \in \Gamma$
be such that the line $pp^*$ contains
the center of the circle $\Gamma$ containing
$\gamma^*$.  Since
$\gamma^*$ is nice, Lemma \ref{circle} 
says that
$\|p-p^*\| \leq \chi(D,d_1) 
\leq \chi(1,d_1).$
\newline

To finish the proof,
we just have to show that $(e,Q)$ is normal.
We enlarge the set
of possible pairs we consider, by allowing
rectangles in $[-3/2,3/2]^2$ having sides parallel
to the coordinate axes and maximum side
length $1/2$.  The same arguments as above,
Lemma \ref{diam} and the $2$-Lipschitz nature
of $\Sigma^{-1}$, show that
$\alpha_1'$ and $\alpha_2'$ are still nice
for any such pair $(e',Q')$.

If $e'$ is the long side of
a $1/2 \times 10^{-100}$ rectangle $\langle Q' \rangle$ contained
in the $10^{-100}$-neighborhood of the coordinate
axes, then $(e',Q')$ is normal:
The arc $\alpha'_1$ is very nearly the
arc of a great circle and the angle between
$\alpha_1'$ and $\alpha_2'$ is very small,
so all arcs ${\cal F\/}'$ 
are all nearly arcs of great circles.
If some choice $(e,Q)$ is not normal, then
by continuity, there is a choice $(e'',Q'')$
in which $\gamma''$ is the diameter of
one of the two boundary arcs of ${\cal F\/}''$.
There is no other way to switch from normal
to not normal.
But this is absurd because the
boundary arcs, $\alpha_1''$ and $\alpha_2''$, are
nice.
\endproof

\subsection{Dot Product Estimates}
\label{sep}

Let $Q$ be a dyadic segment or a good dyadic square.
Let $\delta$ be the hull separation constant of $Q$.
Let $\{q_i\}$ be the points of $Q$.
We make all the same definitions for a second
dyadic square $Q'$.
We define

\begin{equation}
(Q \cdot Q')_{\rm max\/}=\max_{i,j}(\widehat q_i \cdot \widehat q_j')
+\delta+\delta'+\delta\delta'.
\end{equation}
\begin{equation}
(Q \cdot \{\infty\})_{\rm max\/}=\max_{i}\ \widehat q_i \cdot (0,0,1)
\end{equation}

\noindent
{\bf Connectors:\/}
We say that a {\it connector\/} is a line segment connecting
a point on $\widehat{\langle Q \rangle}$ to any of its closest
points in $\langle \widehat Q \rangle$.
We let $\Omega(Q)$ denote the set of connectors
defined relative to $Q$.  By the Hull Approximation Lemma,
each $V \in \Omega(Q)$ has the 
form $W+\delta U$ where $W \in \langle \widehat Q \rangle$
and $\|U\| \leq 1$.  

\begin{lemma}
\label{dotmin}
$V \cdot V' \leq
(Q \cdot Q')_{\rm max\/}$
for all $(V,V') \in \Omega(Q) \times \Omega(Q')$
\end{lemma}

\startproof
Suppose $V \in  \langle \widehat Q \rangle$ and
$V' \in  \langle \widehat Q' \rangle$.
Since the dot product is bilinear, the restriction of the dot product
to the convex polyhedral set
$\langle \widehat Q \rangle \times \langle \widehat Q' \rangle$
takes on its extrema at vertices. Hence
$V \cdot V' \leq
\max_{i,j} q_i \cdot q_i'$.
In this case, we get the desired inequality
whether or not $Q'=\{\infty\}$.

Suppose $Q' \not = \{\infty\}$ and 
$V, V'$ are arbitrary.
We use the decomposition mentioned above:
\begin{equation}
V=W+\delta U, \hskip 30 pt
V'=W'+\delta' U, \hskip 30 pt
W \in \langle \widehat Q \rangle, \hskip 30 pt
W' \in \langle \widehat Q' \rangle.
\end{equation}
But then, by the Cauchy-Schwarz inequality,
$$
|(V \cdot V')-(W \cdot W')|=
|V \cdot \delta' U' + V' \cdot \delta U + \delta U \cdot \delta' U'|
\leq \delta + \delta' + \delta \delta'.$$
The lemma now follows immediately from this
equation, the previous case applied to $W,W'$, and
the triangle inequality.

Suppose that $Q'=\{\infty\}$.
We already know the result when
$V \in \langle \widehat Q \rangle$.
When
$V \in \widehat{\langle Q \rangle}$
we get the better bound above from
Lemma \ref{mono}
and from the fact that
the dot product $V \cdot (0,0,1)$ varies
monotonically with the distance from $V$ to $(0,0,1)$
and {\it vice versa\/}.  Now we know the
result whenever $V$ is the endpoint of a
connector.  By the linearity of the
dot product, the result holds also
when $V$ is an interior point of
a connector.
\endproof

\newpage

%% file: 4energy.tex
\section{The Energy Theorem}

\subsection{Main Result}
\label{ee}
\label{subdivision}

We think of the energy potential
$G=G_k$ as being a function on $(\R^2 \times \infty)^2$,
via the identification $p \leftrightarrow \widehat p$.
Our Energy Theorem below works for any real $k \geq 2$.

Let $\cal Q$ denote the set of dyadic squares
$[-2,2]^2$ together with the dyadic segments in $[0,4]$,
together with $\{\infty\}$.
When $Q=\{\infty\}$ the constants associated to $Q$,
namely the hull separation constant and the convex
hull diameter, are $0$.

Now we are going to define a function
$\epsilon: {\cal Q\/} \times {\cal Q\/} \to [0,\infty)$.
First of all, for notational convenience we set
$\epsilon(Q,Q)=0$ for all $Q$,
When $Q,Q' \in \cal Q$ are unequal, we define
\begin{equation}
\label{EPSILON}
\epsilon(Q,Q')=
\frac{1}{2} k(k-1) T^{k-2}d^2+
2kT^{k-1} \delta
\end{equation}
Here
\begin{itemize}
\item $d$ is the diameter of $\widehat Q$. 
\item $\delta=\delta(Q)$ is the hull approximation constant for $Q$.
See \S \ref{hullsep}.
\item
$T=T(Q,Q')=2+2(Q \cdot Q')_{\rm max\/}$.  See \S \ref{sep}.
\end{itemize}
This is a rational function in the coordinates of $Q$ and $Q'$.
The quantities $d^2$ and $\delta$ are essentially quadratic in
the side-lengths of $Q$ and $Q'$.  Note that
$\epsilon(\{\infty\},Q')=0$ but
$\epsilon(Q,\{\infty\})$ is nonzero when $Q \not = \{\infty\}$.

Let $B=(Q_0,Q_1,Q_2,Q_3)$. For notational
convenience we set $Q_4=\{\infty\}$.  
We define
\begin{equation}
{\bf ERR\/}(B)=\sum_{i=0}^3
\sum_{j=0}^4 \epsilon(Q_i,Q_j).
\end{equation}

\begin{theorem}[Energy]
\label{ENERGY}
$$\min_{v \in \langle B \rangle} {\cal E\/}_k(v) \geq
\min_{v \in B} {\cal E\/}_k(v)-
{\bf ERR\/}(B).$$
\end{theorem}

\begin{corollary}
\label{SUPER}
Suppose that $B$ is a block such that
\begin{equation}
\min_{v \in B} {\cal E\/}(v)-{\bf ERR\/}(B)>{\cal E\/}_k({\rm TBP\/}).
\end{equation}
Then all configurations in $B$ have higher
energy than the TBP.
\end{corollary}

\subsection{The Subdivision Recommendation}

We can write
\begin{equation}
{\bf ERR\/}(B)=\sum_{i=0}^3 {\bf ERR\/}_i(B), \hskip 30 pt
{\bf ERR\/}_i(B)=
\sum_{j=0}^4 \epsilon(Q_i,Q_j).
\end{equation}
We define the
{\it subdivision recommendation\/} 
to be the index
$i \in \{0,1,2,3\}$ for which ${\bf ERR\/}_i(B)$
is maximal. In the extremely unlike event that
two of these terms coincide, we pick the smaller
of the two indices to break the tie.  The subdivision
recommendation feeds into the algorithm
described in \S \ref{ALG}.

The rest of the chapter is devoted to proving
Theorem \ref{ENERGY}.  The next chapter explains
how we use Theorem \ref{ENERGY} in our proof.

\subsection{A Polynomial Inequality}

Theorem \ref{ENERGY} derives from
the case $M=4$ of the following inequality.

\begin{lemma}
\label{ineq}
Let $M \geq 2$ and $k \geq 2$.
Suppose \begin{itemize}
\item $0 \leq x_1 \leq ... \leq x_M$.
\item $\sum_{i=1}^M \lambda_i=1$ and $\lambda_i \geq 0$ for all $i$.
\end{itemize}
Then
\begin{equation}
\label{INEQ}
\sum_{i=1}^M \lambda_i x_i^k-
\bigg( \sum_{I=1}^M \lambda_i x_i\bigg)^k \leq 
\frac{1}{8} k(k-1) x_M^{k-2}\ (x_M-x_1)^2.
\end{equation}
\end{lemma} 

I discovered Lemma \ref{ineq} experimentally

\begin{lemma}
The case $M=2$ of Lemma \ref{ineq} implies the
rest.
\end{lemma}

\startproof
Suppose that $M \geq 3$.
We have one degree of freedom when we keep 
$\sum \lambda_i x_i$ constant and try to vary
$\{\lambda_j\}$ so as to maximize the left
hand side of the inequality.  The right hand
side does not change when we do this, and 
the left hand side varies linearly. Hence,
the left hand size is maximized when
$\lambda_i=0$ for some $i$. But then
any counterexample
to the lemma for $M \geq 3$ gives rise to
a counter example for $M-1$.
\endproof

In the case $M=2$, we set
$a=\lambda_1$.
Both sides of the inequality in
Lemma \ref{ineq} are homogeneous of
degree $k$, so it suffices to consider the
case when $x_2=1$.  We set $x=x_1$.
The inequality of is then
$f(x) \leq g(x)$, where
\begin{equation}
\label{McM}
f(x)=
(a x^k + 1-a)-(a x + 1-a)^k;
\hskip 30 pt
g(x)=\frac{1}{8}k(k-1)(1-x)^2.
\end{equation}
This is supposed to hold for
all $a, x \in [0,1]$.

The following argument is due to
C. McMullen, who figured it out
after I told him about the inequality.

\begin{lemma}
Equation \ref{McM} holds for
all $a,x \in [0,1]$ and all $k \geq 2$.
\end{lemma}

\startproof
Equation \ref{McM}. We think of $f$ as
a function of $x$, with $a$ held fixed.
Since $f(1)=g(1)=1$, it suffices
to prove that $f'(x) \geq g'(x)$ on $[0,1]$.
Define
\begin{equation}
\phi(x)=akx^{k-1}, \hskip 30 pt
b=(1-a)(1-x).
\end{equation}
We have
\begin{equation}
-f'(x)=\phi(x+b) - \phi(x).
\end{equation}
Both $x$ and $x+b$ lie in $[0,1]$. So,
by the mean value theorem there is some
$y \in [0,1]$ so that
\begin{equation}
\frac{\phi(x+b)-\phi(x)}{b}=\phi'(y)=
ak(k-1)y^{k-2}.
\end{equation}
Hence
\begin{equation}
-f'(x)=b\phi'(y)=a(1-a)k(k-1)(1-x)y^{k-2}
\end{equation}
But
$a(1-a) \in [0,1/4]$ and
$y^{k-2} \in [0,1].$
Hence
\begin{equation}
-f'(x) \leq \frac{1}{4}k(k-1)(1-x)=-g'(x).
\end{equation}
Hence $f'(x) \geq g'(x)$ for all $x \in [0,1]$.
\endproof

\noindent
{\bf Remark:\/} Lemma \ref{ineq} has the following
motivation.  The idea behind the Energy Theorem is that
we want to measure the deviation of the energy function
from being linear, and for this we would like a
quadratic estimate.  Since our energy $G_k$ involves
high powers, we want to estimate these high powers by
quadratic terms.

\subsection{The Local Energy Lemma}

Let $Q=\{q_1,q_2,q_3,q_4\}$ be the vertex set 
of $Q \in \cal Q$. We allow
for the degenerate case that $Q$ is a
line segment or $\{\infty\}$.  In this case we
just list the vertices multiple times,
for notational convenience.

Note that every point in the convex
quadrilateral $\langle \widehat Q \rangle$
is a convex average of the vertices.
For each $z \in \langle Q \rangle$, there is a
some point $z^* \in \langle \widehat Q\rangle$
which is as close as possible to
$\widehat z \in \widehat{\langle Q \rangle}$.  There are
constants $\lambda_i(z)$ such that
\begin{equation}
\label{zstar}
\label{starz}
z^*=\sum_{i=1}^4 \lambda_i(z)\ \widehat q_j,
\hskip 30 pt
\sum_{i=1}^4 \lambda_i(z)=1.
\end{equation}
We think of the $4$ functions
$\{\lambda_i\}$ as a partition of unity on $\langle Q \rangle$.
The choices above might not be unique,
but we make such choices once and for all for each $Q$.
We call the assignment $Q \to \{\lambda_i\}$
the {\it stereographic weighting system\/}.

\begin{lemma}[Local Energy]
Let $\epsilon$ be the function defined in
the Energy Theorem.
Let $Q,Q'$ be distinct members of $\cal Q$. Given
any $z \in Q$ and $z' \in Q'$, 
\begin{equation}
G(z,z') \geq \sum_{i=1}^4 \lambda_i(z)G(q_i,z')-\epsilon(Q,Q').
\end{equation}
\end{lemma}

\startproof
For notational convenience, we set $w=z'$.
Let
\begin{equation}
X=\big(2+2z^*\cdot \widehat w\big)^k.
\end{equation}
The Local Energy Lemma follows from adding
these two inequalities:
\begin{equation}
\label{zoop}
\sum_{i=1}^4 \lambda_i G(q_i,w) - X \leq
\frac{1}{2} k(k-1)T^{k-2}d^2
\end{equation}
\begin{equation}
\label{zoop2}
X-G(z,w) \leq
2kT^{k-1}\delta.
\end{equation}
We will establish these inequalities in turn.
\newline
\newline
Let $q_1,q_2,q_3,q_4$ be the vertices of $Q$.
Let $\lambda_i=\lambda_i(z)$.  
We set
\begin{equation}
x_i=4-\|\widehat q_i-\widehat w\|^2=2+2\widehat q_i \cdot \widehat w, \hskip 30 pt i=1,2,3,4.
\end{equation}
Note that $x_i \geq 0$ for all $i$. 
We order so that
$x_1 \leq x_2 \leq x_3 \leq x_4$. We have
\begin{equation}
\label{term1}
\sum_{i=1}^4 \lambda_i(z)G(q_i,w)=\sum_{i=1}^4 \lambda_i x_i^k,
\end{equation}
\begin{equation}
\label{term2}
X=(2+2\widehat z^* \cdot \widehat w)^k=
\bigg(\sum_{i=1}^4 \lambda_i (2+\widehat q_i \cdot \widehat w)\bigg)^k=
\bigg(\sum_{i=1} \lambda_i x_i\bigg)^k.
\end{equation}

Combining Equation \ref{term1}, Equation \ref{term2}, and
the case $M=4$ of Lemma \ref{ineq},
\begin{equation}
\label{boundX}
\sum_{i=1}^4 \lambda_i G(q_i,w) - X=
\sum_{i=1}^4 \lambda_i x_i^k - \bigg(\sum_{i=1}^4 \lambda_i x_i\bigg)^k
\leq
\frac{1}{8} k(k-1)x_4^{k-2}(x_4-x_1)^2.
\end{equation}

By Lemma \ref{dotmin}, we have
\begin{equation}
\label{subs1}
x_4=2+2(\widehat q_4 \cdot \widehat w) \leq
2+2(Q \cdot Q')_{\rm max\/}=T.
\end{equation}
Since $d$ is the diameter of  $\langle \widehat Q \rangle$
and $\widehat w$ is a unit vector, 
\begin{equation}
\label{subs2}
x_4-x_1=2 \widehat w \cdot (\widehat q_4-\widehat q_1)
\leq 2\|\widehat w\| \|\widehat q_4-\widehat q_1\|=
2 \|\widehat q_4-\widehat q_1\| \leq 2d.
\end{equation}
Plugging Equations \ref{subs1} and
\ref{subs2} into Equation
\ref{boundX}, we get Equation \ref{zoop}.

Now we establish Equation \ref{zoop2}.
Let $\gamma$ denote the unit speed line segment connecting
$\widehat z$ to $z^*$.  Note that the length $L$ of
$\gamma$ is at most $\delta$, by the Hull Approximation Lemma.
Define
\begin{equation}
f(t)=
\bigg(2+2 \widehat w \cdot \gamma(t)\bigg)^k.
\end{equation}
We have $f(0)=X$. Since $\widehat w$ and $\gamma(1)=\widehat z$
are unit vectors, $f(L)=G(z,w)$. Hence
\begin{equation}
\label{integ}
X-G(z,w)=f(0)-f(L), \hskip 30 pt L \leq \delta.
\end{equation}
By the Chain Rule,
\begin{equation}
f'(t)=\big(2 \widehat w \cdot \gamma'(t)\big) \times
k\bigg(2+2\widehat w \cdot \gamma(t)\bigg)^{k-1}.
\end{equation}
Note that $|2\widehat w \cdot \gamma'(t)| \leq 2$
because both vectors are unit vectors.
Note that $\gamma$ parametrizes one of the connectors from
Lemma \ref{dotmin}, so we have
\begin{equation}
\label{diff}
|f'(t)| \leq 2k\bigg(2+2\widehat w \cdot \gamma(t)\bigg)^{k-1}  \leq 
2k\bigg(2+2(Q \cdot Q')_{\rm max\/}\bigg)^{k-1}=2kT^{k-1}.
\end{equation}
Equation \ref{zoop2} now follows from Equation \ref{integ},
Equation \ref{diff}, and integration.
\endproof

\subsection{From Local to Global}

Let $\epsilon$ be the function from
the Energy Theorem.  
Let $B=(Q_0,...,Q_N)$ be a list of
$N+1$ elements of $\cal Q$.  We care about
the case $N=4$ and $Q_4=\{\infty\}$, but
the added generality makes things clearer.
Let 
$q_{i,1},q_{i,2},q_{i,3},q_{i,4}$ be the vertices of $Q_i$.
The vertices of $\langle B \rangle$ are indexed by
a multi-index $$I=(i_0,...,i_n) \in \{1,2,3,4\}^{N+1}.$$
Given such a multi-index, which amounts to a choice
of vertex of $\langle B \rangle$, 
we define the energy of the corresponding vertex configuration:
\begin{equation}
{\cal E\/}(I)=
{\cal E\/}(q_{0,i_0},...,q_{N,i_N})
\end{equation}
We will prove the following sharper result.

\begin{theorem}
\label{ENERGY2}
Let $z_0,...,z_N \in \langle B \rangle$.
Then
\begin{equation}
\label{AVE}
{\cal E\/}(z_0,...,z_N) \geq
\sum_{I} \lambda_{i_0}(z_0)...\lambda_{i_N}(z_N){\cal E\/}(I)-
\sum_{i=0}^N \sum_{j=0}^N \epsilon(Q_i,Q_j).
\end{equation}
The sum is taken over all multi-indices.
\end{theorem}

\begin{lemma}
Theorem \ref{ENERGY2} implies the
Energy Theorem.
\end{lemma}

\startproof
Notice that
\begin{equation}
\sum_{I} \lambda_{i_0}(z_0)...\lambda_{i_N}(z_N)=
\prod_{j=0}^N \bigg(\sum_{a=1}^4 \lambda_a(z_j)\bigg)=1.
\end{equation}
Therefore
\begin{equation}
\sum_{I} \lambda_{i_0}(z_0)...\lambda_{i_N}(z_N){\cal E\/}(I)
\geq \min_{v \in B} {\cal E\/}(v),
\end{equation}
because the sum on the left hand side is the convex 
average of vertex energies and the term on the right
is the minimum of the vertex energies.

For any $(z_0,...,z_N) \in \langle B \rangle$, 
we now know from Theorem \ref{ENERGY2} that
$$
{\cal E\/}(z_0,...,z_N) \geq
 \min_{v \in B} {\cal E\/}(v)-
\sum_{i=0}^N \sum_{j=0}^N \epsilon(Q_i,Q_j).
$$
When we take $N=4$ and $Q_4=\{\infty\}$, the
second expression on the right hand side of
this last equation is precisely
${\bf ERR\/}(B)$.  This establishes
the Energy Theorem.
\endproof

We now prove Theorem \ref{ENERGY2}.

\subsubsection{A Warmup Case}

Consider the case when $N=1$. 
Setting $\epsilon_{ij}=\epsilon(Q_i,Q_j)$, the
Local Energy Lemma gives us
\begin{equation}
G(z_0,z_1) \geq 
\sum_{\alpha=1}^4 \lambda_{\alpha}(z_0)G(q_{0\alpha},z_1)-\epsilon_{01}.
\end{equation}
\begin{equation}
G(q_{0\alpha},z_1) \geq
\sum_{\beta=1}^4 \lambda_{\beta}(z_1)G(q_{1\beta}(z_1),q_{0\alpha}) - \epsilon_{10}.
\end{equation}
Plugging the second equation into the first and using 
$\sum \lambda_{\alpha}(z_0)=1$, we have
\begin{equation}
\label{AVE1}
G(z_0,z_1) \geq \bigg(\sum_{\alpha,\beta} \lambda_{\alpha}(z_0)\lambda_{\beta}(z_1)
G(q_{0\alpha},q_{1\beta})\bigg) - (\epsilon_{01}+\epsilon_{10}).
\end{equation}
This is precisely Equation \ref{AVE} when $N=1$.

\subsubsection{The General Case}

Now assume that $N \geq 2$.
We rewrite Equation \ref{AVE1} as follows:
\begin{equation}
\label{AVE2}
G(z_0,z_1) \geq \sum_{A} \lambda_{A_0}(z_0)\lambda_{A_1}(z_1)\ G(q_{0A_0},q_{1A_1})-
(\epsilon_{01}+\epsilon_{10}).
\end{equation}
The sum is taken over multi-indices $A$ of length $2$.

We also observe that
\begin{equation}
\label{AVE3}
\sum_{I'} \lambda_{i_2}(z_2)...\lambda_{i_N}(z_N)=1.
\end{equation}
The sum is taken over all multi-indices
$I'=(i_2,...,i_N)$.
Therefore, if $A$ is held fixed, we have
\begin{equation}
\lambda_{A_0}(z_0)\lambda_{A_1}(z_1) =
\sum'_{I} \lambda_{I_0}(z_0)...\lambda_{I_N}(z_N).
\end{equation}
The sum is taken over all multi-indices of length $N+1$ which 
have $I_0=A_0$ and $I_1=A_1$.
Combining these equations, we have
\begin{equation}
G(z_0,z_1) \geq 
\sum_I \lambda_{I_0}(z_0)... \lambda_{I_N}(z_N) G(q_{0I_0},q_{1I_1})-
(\epsilon_{01}+\epsilon_{10}).
\end{equation}
The same argument works for other pairs of indices, giving
\begin{equation}
\label{AVE4}
G(z_i,z_j) \geq 
\sum_I \lambda_{I_0}(z_0)... \lambda_{I_N}(z_N) G(q_{iI_i},q_{jI_j})-
(\epsilon_{ij}+\epsilon_{ji}).
\end{equation}

Now we interchange the order of summation and observe that
$$
\sum_{i<j}
\bigg( \sum_I \lambda_{I_0}(z_0)...\lambda_{I_N}(z_N)\ G(q_{iI_i},q_{jI_j})\bigg)=$$
$$
\sum_I \sum_{i<j} \lambda_{I_0}(z_0)... \lambda_{I_N}(z_N)\ G(q_{iI_i},q_{jI_j})=
$$
$$
\sum_I \lambda_{I_0}(z_0)... \lambda_{I_N}(z_N) 
\Bigg(\sum_{i<j} G(q_{iI_i},q_{jI_j})\Bigg)=
$$
\begin{equation}
\sum_I \lambda_{I_0}(z_0)... \lambda_{I_N}(z_N)\ {\cal E\/}(I).
\end{equation}

Therefore, when we sum Equation \ref{AVE4} over all $i<j$, we
get precisely the inequality in 
Equation \ref{AVE}.  This completes the proof.

\subsection{A More General Result}
\label{energyX}

Though we have no need for it in this paper, we
mention an easy generalization of the Energy Theorem.
Suppose we have some energy of the form
\begin{equation}
F=\sum_{k=1}^N a_k G_k
\end{equation}
where $a_1,...,a_N$ is some sequence of
real numbers, not necessarily positive.
We define
\begin{equation}
{\bf ERR\/}_F=\sum_{k=1}^N |a_k|\ {\bf ERR\/}_k
\end{equation}
Here ${\bf ERR\/}_k$ is the error term associated
to $G_k$ in the Energy Theorem.
With this definition in place, we have

\begin{theorem}
\label{ENERGYxx}
$$\min_{v \in \langle B \rangle} {\cal F\/}_k(v) \geq
\min_{v \in B} {\cal F\/}_k(v)-
{\bf ERR\/}_F(B).$$
\end{theorem}

\newpage

%% file: 5divide.tex
\section{The Algorithm}
\label{ALG}

\subsection{Grading a Block}
\label{GRADE}

In this section we describe what we mean
by {\it grading\/} a block.  This step
feeds into the computational algorithm
presented in the next section.

Let $B_0$ denote the cube of in-radius
$\epsilon_0=2^{-15}$ about the
configuration of $\square$ representing
the normalized TBP.
We fix some energy $G_k$.  We perform the
following tests on a block 
$B=(Q_0,Q_1,Q_2,Q_3)$, in the order listed.
\begin{enumerate}

\item If the calculations in Equation \ref{irrelevantX}
deem $B$ irrelevant, we pass $B$.

\item If some component square $Q_i$ of $B$
has side length more than $1/2$ we fail $B$ and
recommend that $B$ be subdivided along the
first such index.

\item If we compute that $Q_i \not \subset
[-3/2,3/2]^2$ for some $i=1,2,3$, we pass $B$.
Given the previous step, $Q_i$ is disjoint
from $(-3/2,3/2)^2$.

\item If the calculations in \S \ref{checkTBP} show
that $B \subset B_0$, we pass $B$.  Here we
take $S=2^{25}$ and (as we have already said)
$\epsilon_0=2^{-15}$.

\item If the calculations in
\S \ref{ee} show that
$B$ satisfies Corollary \ref{SUPER}, we
pass $B$.
Otherwise, we fail $B$ and pass along the
recommended subdivision.
\end{enumerate}

\subsection{Depth First Search}

We plug the grading step into the following
algorithm.

\begin{enumerate}
\item Begin with a list LIST of blocks in $\square$.
Initially LIST consists of a single element,
namely $\square$.

\item Let $B$ be the last member of LIST. We delete
$B$ from LIST and then we grade $B$.

\item Suppose $B$ passes.  If
LIST is empty, we halt and declare success.
Otherwise, we return to Step 2.

\item Suppose $B$ fails.  
In this case, we subdivide $B$ along the
subdivision recommendation and we
append to LIST the
subdivision of $B$.  Then we return
to Step 2.
\end{enumerate}

If the algorithm halts with success, it implies
that every relevant block $B$ either lies
in $B_0$ or does not contain a minimizer.

\subsection{The Results}

I will detail the
technical implementations of the
algorithm in \S \ref{computer}.

For each $k=3,4,5,6$, 
I ran the programs (most recently) on August 4-5, 2016
on my 2014 Macbook pro.  
\begin{itemize}

\item For $G_3$ the program finished in about
$1$ hour and $25$ minutes

\item For $G_4$ the program finished in about
$1$ hour and $39$ minutes

\item For $G_5$ the program finished in about
$2$ hours and $40$ minutes.

\item For $G_6$ the program finished in about
$5$ hours and $38$ minutes.
\end{itemize}

In each case, the program
produces a partition of $\square$
into $N_k$ smaller blocks, each of
which is either irrelevant, contains no
minimizer, or lies in $B_0$. Here
$$
(N_3,N_4,N_5,N_6)=
(5513537,6201133,9771906,20854602).
$$

These calculations rigorously establish the
following result.

\begin{lemma}
\label{energy1}
Let $B_0 \subset \square$ denote the cube of 
in-radius $2^{-15}$ about $P_0$.
If $P \in \square$ is a minimizer with
respect to any of
$G_3,G_4,G_5,G_6$ then $P$ has
the same energy as a
configuration in $B_0$.
\end{lemma}

The algorithm runs for $G^{\#}_{10}$
with the following modifications.
\begin{itemize}
\item We use $\epsilon_0=2^{-18}$ and
$S=2^{30}$.  
\item We use Theorem \ref{ENERGYxx} in place
of Theorem \ref{ENERGY}.
\end{itemize}

I ran the algorithm on $G^{\#}_{10}$ in
the last week of July 2016 on my $2014$ iMac.
The calculation ran to completion after
about $51$ hours and $13$ minutes, producing
a partition of size
$67899862$.   The calculation establishes
the following result.

\begin{lemma}
\label{energy2}
Let $B_0^{\#} \subset \square$ denote the cube of 
in-radius $2^{-18}$ about $P_0$.
If $P \in \square$ is a minimizer with
respect to $G^{\#}_{10}$ then
$P$ has the same energy as a
configuration in $B_0^{\#}$.
\end{lemma}

\newpage

%% file: 6hessian.tex
\section{Local Analysis of the Hessian}
\label{local}

\subsection{Eigenvalues of Symmetric Matrices}
\label{eig00}

Let $H$ be a symmetric $n \times n$ real matrix.
$H$ always has an orthonormal basis of eigenvectors,
and real eigenvalues.  $H$ is {\it positive definite\/}
if all these eigenvalues are positive.  This is equivalent
to the condition that $Hv \cdot v>0$ for all nonzero $v$.
More generally, $Hv \cdot v \geq \lambda \|v\|$, where
$\lambda$ is the lowest eigenvalue of $H$.  

Here is one way to bound the lowest eigenvalue
of $H$.

\begin{lemma}[Alternating Criterion]
Let $\chi(t)$ be the characteristic polynomial of $H$.
Suppose that the coefficients of $P(t)=\chi(t+\lambda)$ are
alternating and nontrivial.  Then the lowest eigenvalue
of $H$ exceeds $\lambda$.
\end{lemma}

\startproof
An alternating polynomial has no negative roots.
So, if $\chi(t+\lambda)=0$ then $t>0$ and
$t+\lambda>\lambda$.
\endproof

Let $H_0$ be some positive definite symmetric
matrix and let $\Delta$ be
some other symmetric matrix of the same size.
Recall various definitions of the $L_2$ matrix norm:
\begin{equation}
\|\Delta\|_2=\sqrt{\sum_{ij}\Delta_{ij}^2}=
\sqrt{{\rm Trace\/}(\Delta\Delta^t)}=
\sup_{\|v\|=1} \|\Delta v\|.
\end{equation}

\begin{lemma}[Variation Criterion]
Suppose that $\|\Delta\|_2 \leq \lambda$, where
$\lambda$ is some number less than the lowest
eigenvalue
of $H_0$.  Then $H=H_0+\Delta$ is also positive definite.
\end{lemma}

\startproof
$H$ is positive definite if and only of
$Hv \cdot v>0$ for every nonzero unit vector $v$.
Let $v$ be such a vector. Writing $v$ in an
orthonormal basis of eigenvectors we see that
$H_0v \cdot v>\lambda$. Hence
$$
Hv \cdot v = (H_0 v+\Delta v) \cdot v  \geq
H_0v \cdot v - |\Delta v \cdot v|>
\lambda - \|\Delta v\| \geq
\lambda - \|\Delta\|_2 \geq 0.
$$
This completes the proof.
\endproof

\subsection{Taylor's Theorem with Remainder}
\label{taylor1}

In this section we are just packaging a special case of
Taylor's Theorem with Remainder.  Here are 
some preliminary definitions.

\begin{itemize}
\item Let $P_0 \in \R^7$ be some point.
\item Let $B$ denote some cube of in-radius
$\epsilon$ centered at $P_0$.
\item $\phi: \R^7 \to \R$ be some function.
\item Let
$\partial_I \phi$ be the
partial derivative of $\phi$ w.r.t. a
multi-index $I=(i_1,...,i_7)$. 
\item  Let
$|I|=i_1+\cdots i_7$.  This is the {\it weight\/} of $I$.
\item Let $I!=i_1! \cdots i_7!$.
\item Let $\Delta^I=x^{i_1}...x^{i_7}$.
Here $\Delta=(x_1,...,x_7)$ is some vector.
\item For each positive integer $N$ let
\begin{equation}
\label{biggie}
M_N(\phi)=\sup_{|I|=N}\sup_{P \in B} |\partial_I \phi(P)|,
\hskip 30 pt
\mu_N(\phi)=\sup_{|I|=N}|\partial_I \phi(P_0)|.
\end{equation}
\end{itemize}

Let $U$ be some open neighborhood of $B$.
Given $P \in B$, let $\Delta=P-P_0$.
Taylor's Theorem with Remainder says that there
is some $c \in (0,1)$ such that
$$
\phi(P)= \sum_{a=0}^N 
\sum_{|I|=a} \frac{|\partial_I \phi(P_0)|}{I!}\Delta^I+
\sum_{|I|=N+1} \frac{\partial_I f(P_0+c\Delta)}{I!}\Delta^I
$$

Using the fact that
$$
|\Delta^I| \leq \epsilon^{|I|}, \hskip 30 pt
\sum_{|I|=m} \frac{1}{I!}=\frac{7^{m}}{m!},
$$
and setting $N=4$ we get
\begin{equation}
\label{taylor}
\sup_{P \in B} |\phi(P)| \leq |\phi(P_0)|+
\sum_{j=1}^4 \frac{(7\epsilon)^j}{j!} \mu_j(\phi) +
\frac{(7\epsilon)^5}{(5)!}M_5(\phi).
\end{equation}

\subsection{The Lowest Eigenvalues}
\label{TBP2}

We will first deal with $G_k$ for $k=3,4,5,6$ and
then at the end of the chapter we will explain
the modifications needed to handle $G^{\#}_{10}$.
Let 
\begin{equation}
\epsilon_0=2^{-15}
\end{equation}
 and $B_0$ be
as in Lemma \ref{energy1}
Unless otherwise stated, we will
take $k \in \{2,3,4,5,6\}$.
(We treat the case $k=2$ because it will
be useful to us when we need to
deal with $G^{\#}_{10}$.)
We have the energy map
$E_k: B_0 \to \R_+$
given by
\begin{equation}
E_k(x_1,...,x_7)=
\sum_{i<j} G_k(\Sigma^{-1}(p_i)-\Sigma^{-1}(p_j)).
\end{equation}
Here we have set $p_4=\infty$, and
$p_0=(x_1,0)$ and $p_i=(x_{2i},x_{2i+1})$ for $i=1,2,3$.
As usual $\Sigma$ is stereographic projection.

Let $P_0 \in B_0$ denote the point
corresponding to the TBP.  It follows
from symmetry, and also from a direct
calculation, that $P_0$ is a critical
point for $E_k$ in all cases.
The Main Theorem for $G_k$ now follows
from Lemma \ref{energy1} and from
the statement that the Hessian of $G_k$
is positive definite throughout $B_0$.

Let $H_k$ denote the Hessian of $E_k$.
Define
\begin{equation}
(\lambda_2,\lambda_3,\lambda_4,\lambda_5,\lambda_6)=
(1,10,24,36,43).
\end{equation}

\begin{lemma}
\label{base}
For each $k=2,3,4,5,6$, the Hessian $H_k(P_0)$ is positive
definite and its lowest eigenvalue exceeds $\lambda_k$.
\end{lemma}

\startproof
All the pairs $(H_k,\lambda_k)$ satisfy
the Alternating Criterion from \S \ref{eig00}.
\endproof

\subsection{The Target Bounds}

Now we set up the inequalities we need
to establish in order to show that
the hessians $H_k$ are positive
definite throughout the small neighborhood $B_0$.

Define
\begin{equation}
\label{thirdbound}
F_k=\sqrt{\sum_{|J|=3} M_{J,k}^2}, \hskip 30 pt
M_{J,k}=\sup_{P \in B_0} |\partial_J E_k(P)|.
\end{equation}
The sum is taken over all multi-indices $J$ of
weight $3$.

\begin{lemma}
\label{smallL2}
Suppose that
$\sqrt 7 \epsilon_0 F_k \leq \lambda_k$.
Then the Hessian of $E_k$ is positive
definite throughout the neighborhood $B_0$.
\end{lemma}

\startproof
We suppress the value of $k$ from our notation.
Let $H_0$ denote the Hessian of $E$ at $P_0$.
Let $H$ denote the Hessian of $E$ at $P$.
Let $\Delta=H-H_0$.  Clearly
$H=H_0+\Delta$.

Let $\gamma$ be the unit
speed line segment connecting $P$ to $P_0$ in $\R^7$.
Note that $\gamma \subset B_0$ and
$\gamma$ has length 
$L \leq \sqrt 7 \epsilon_0$.
We set $H_{L}=H$ and we let
$H_t$ be the Hessian of $E$ at the point of
$\gamma$ that is $t$ units from $H_0$.

We have
\begin{equation}
\Delta=\int_0^{L} D_t(H_t)\ dt.
\end{equation}
Here $D_t$ is the unit directional
derivative of $H_t$ along $\gamma$.

Let $(H_t)_{ij}$ denote the $ij$th entry of
$H_t$.  Let $(\gamma_1,...,\gamma_7)$ be the
components of the unit vector in the direction of $\gamma$.
Using the fact that $\sum_k \gamma_k^2=1$ and
the Cauchy-Schwarz inequality, and the
fact that mixed partials commute, we have
\begin{equation}
(D_tH_t)_{ij}^2=
\bigg(\sum_{k=1}^7 \gamma_k 
\frac{\partial}{\partial x_k}
\frac{\partial^2 H_t}{\partial x_i \partial x_j}\bigg)^2 \leq
\sum_{k=1}^7 
\bigg(\frac{\partial^3 H_t}{\partial x_i \partial x_j \partial x_k}\bigg)^2.
\end{equation}
Summing this inequality over $i$ and $j$ we get
\begin{equation}
\|D_tH_t\|_2^2 \leq
\sum_{i,j,k}
\bigg(\frac{\partial^3 H_t}{\partial x_i \partial x_j \partial x_k}\bigg)^2 \leq F^2.
\end{equation}
Hence
\begin{equation}
\label{small2}
\|\Delta\|_2 \leq \int_0^L \|D_t(H_t)\|_2\ dt \leq
L F \leq \sqrt 7 \epsilon_0 F<\lambda.
\end{equation}
This lemma now follows immediately from
Lemma \ref{smallL2}.
\endproof

Referring to Equation \ref{biggie}, and
with respect to the neighborhood $B_0$, define
$M_8(E_k)$ and $\mu_j(E_k)$ for $j=4,5,6,7$. 
Let $J$ be any multi-index of weight $3$. 
Using the fact that
$$
\mu_j(\partial_J E_K) \leq \mu_{j+3}(E_k), \hskip 30 pt
M_{5}(\partial_J E_K) \leq M_{8}(E_k),
$$
we see that Equation \ref{taylor} gives us the
bound
\begin{equation}
\label{BOUND}
M_{J,k} \leq |\partial_JE_k(P_0)| + 
\sum_{j=1}^4
\frac{(7\epsilon_0)^j}{j!}\mu_{j+3}(E_k)
+\frac{(7\epsilon_0)^5}{5!}M_8(E_k).
\end{equation}

\subsection{The Biggest Term}
\label{combinatorics}

In this section we will prove that the last term
in Equation \ref{BOUND} is at most $1$.
Recall that $\Sigma$ is stereographic projection.
Define

\begin{equation}
f_k(a,b)=\bigg(4-\|\Sigma^{-1}(a,b)-(0,0,1)\|^2\bigg)^k=
4^k\bigg(\frac{a^2+b^2}{1+a^2+b^2}\bigg)^k.
\end{equation}

$$
g_k(a,b,c,d)=\bigg(4-\|\Sigma^{-1}(a,b)-\Sigma^{-1}(c,d)\|^2\bigg)^k=$$
\begin{equation}
4^k\bigg(\frac{1 + 2 a c + 2 b d + (a^2+b^2)(c^2+d^2) }
{(1 + a^2 + b^2) (1 + c^2 + d^2)}\bigg)^k
\end{equation}
Note that
\begin{equation}
\label{shortcut}
f_k(a,b)=\lim_{c^2+d^2 \to \infty} g_k(a,b,c,d).
\end{equation}
We have
$$
E_k(x_1,...,x_7)=
f_k(x_1,0)+f_k(x_2,x_3)+f_k(x_4,x_5)+f_k(x_6,x_7)+$$
$$
g_k(x_1,0,x_2,x_3)+
g_k(x_1,0,x_4,x_5)+
g_k(x_1,0,x_6,x_7)+$$
$$
g_k(x_2,x_3,x_4,x_5)+
g_k(x_2,x_3,x_6,x_7)+
g_k(x_4,x_5,x_6,x_7).
$$
Each variable appears in at most $4$
terms, $3$ of which appear in a
$g$-function and $1$ of which appears
in an $f$-function. Hence
\begin{equation}
\label{bound1}
M_8(E_k) \leq M_8(f_k)+3M_8(g_k) \leq 4M_8(g_k).
\end{equation}
The last inequality is a consequence of Equation \ref{shortcut}
and we use it so that we can concentrate on just one
of the two functions above.

\begin{lemma}
\label{easy}
When $r,s,D$ are non-negative integers and $r+s \leq 2D$,
$$\bigg{|}\frac{x^ry^s}{(1+x^2+y^2)^{D}}\bigg{|}<1.$$
\end{lemma}

\startproof
The quantity factors into expressions of the form
$|x^{\alpha}y^{\beta}/(1+x^2+y^2)|$ where
$\alpha+\beta \leq 2$. 
Such quantities are bounded above by $1$.
\endproof

For any polynomial $\Pi$, let $|\Pi|$ denote
the sum of the absolute values of the
coefficience of $\Pi$.
For each $8$th derivative $D_Ig_k$, we have
\begin{equation}
D_Ig_k=\frac{\Pi(a,b,c,d)}{(1+a^2+a^2)^{k+8}(1+c^2+d^2)^{k+8}},
\end{equation}
Where $\Pi_I$ is a polynomial of maximum
$(a,b)$ degree at most $2k+16$ and maximum
$(c,d)$ degree at most $2k+16$.
Lemma \ref{easy} then gives 
\begin{equation}
\sup_{(a,b,c,d) \in \R^4} |D_Ig_k(a,b,c,d)| \leq |\Pi_I|.
\end{equation}
Thus, we compute in Mathematica that
\begin{equation}
\label{bound2}
M_8(g_k) \leq
\sup_{k=2,3,4,5,6} \sup_I |\Pi_I|=13400293856913653760<2^{64}.
\end{equation}
The max is achieved when $k=6$ and $I=(8,0,0,0)$ or
$(0,8,0,0)$, etc.

Combining Equations \ref{bound1} and
\ref{bound2}, we get
\begin{equation}
\frac{\bigg(7\times2^{-15}\bigg)^5}{5!} \times M_8(E_k)<1,
\end{equation}
$k=2,3,4,5,6$.

\subsection{The End of the Proof}
\label{exception}

Let
\begin{equation}
\mu^*_{j+3,k}=\frac{(7\epsilon_0)^{j}}{(j)!} \mu_{j+3}(E_k)
\end{equation}
We now estimate
$\mu^*_{j,k}$ for $j=4,5,6,7$ and
$k=2,3,4,5,6,10$.

The same considerations as above show that
\begin{equation}
\mu^*_{j+3,k} \leq \frac{(7\epsilon_0)^{j}}{(j)!}
\bigg(\mu_{j+3}(f_k)+ 3\mu_{j+3}(g_k)\bigg).
\end{equation}
Here we are evaluating the $(j+3)$rd partial
 at all points which
arise in the TBP configuration and then
taking the maximum.  For instance,
for $g_k$, one choice would be
$(a,b,c,d)=(1,0,-1/2,\sqrt 3/2)$.
Here is a computed matrix of upper bounds for
$\mu^*_{j,k}$.  The rows give the fixed exponent $k$.
\begin{equation}
\label{matxX}
\matrix{j:&4&5&6&7\cr
k=2:&1&1&1&1 \cr
k=3:&3&1&1&1 \cr k=4:&9 &1&1&1 \cr k=5:&12 &1&1 &1 \cr
k=6:&122&1&1&1 \cr
k=10:& 47480&44&1&1}
\end{equation}
The first and last rows are needed for the next section.

Given these bounds, Equation \ref{BOUND} gives
\begin{equation}
\label{GREAT}
|M_{J,k}|< 1+ \mu_7^*+\mu_6^*+\mu_5^*+\mu_4^*<
1+1+1+1+122<500.
\end{equation}
The first $1$ comes from the previous section.
We have been generous in the last inequality
to illustrate that the situation is not delicate here.

For any real vector $V=(V_1,...,V_{343})$ define
\begin{equation}
\overline V=(|V_1|+500,...,|V_{343}|+500).
\end{equation}

Let $V_k$ denote the vector of
third partials of $E_k$, evaluated at $P_0$, and
ordered (say) lexicographically.
In view of Equation \ref{GREAT},
we have
$F_k \leq \|\overline V_k\|$. 
Hence
\begin{equation}
\sqrt 7\epsilon_0 F_k \leq  \sqrt 7\epsilon_0 \|\overline V_k\|.
\end{equation}
Rounding up to the nearest integer, we compute
the coordinatewise inequality
\begin{equation}
\label{bounds3456}
\sqrt 7\epsilon_0 
(\|\overline V_2\|,\|\overline V_3\|,...,\|\overline V_6\|)<
(1,1,2,4,12) \leq (\lambda_2,\lambda_3,\lambda_4,\lambda_5,\lambda_6).
\end{equation}
This completes the proof of the Main Theorem for $G_3,G_4,G_5,G_6$.

\subsection{The Last Case}

Now we deal with $G^{\#}_{10}$.
Let $\epsilon_0^{\#}=2^{-18}$.
The pair $(H^{\#}_{10},1448)$
satisfies the Alternating Criterion.
Hence, the Hessian $H^{\#}_{10}$
has lowest eigenvalue greater than $1448$.
In view of Lemma \ref{energy2} and
the results in this chapter,
the Main Theorem for $G^{\#}_{10}$ follows
from the inequality
\begin{equation}
\sqrt 7 \epsilon_0^{\#} F_{10} +
28 \sqrt 7 \epsilon_0^{\#} F_{5} +
102 \sqrt 7 \epsilon_0^{\#} F_{2}<1448.
\end{equation}
  Since $\epsilon_0^{\#}<\epsilon_0$,
the bounds in Equation \ref{bounds3456} remains
true with respect to $\epsilon^{\#}_0$.
This reduces our goal to showing that
\begin{equation}
\sqrt 7 \epsilon_0 F_{10}<1234
\hskip 20 pt (=1448-28 \times 4 - 102 \times 1).
\end{equation}

We compute that
\begin{equation}
\label{badbound}
M_8(g_{10}) \leq 162516942801336639946752000<2^{88}
\end{equation}
This leads to
\begin{equation}
\frac{\bigg(7\times 2^{-18}\bigg)^5}{5!} \times M_8(E_{10})<19
\end{equation}

Since we are using $\epsilon^{\#}_0$ in
place of $\epsilon_0$, we can divide the
last row of the matrix in Equation \ref{matxX} by
$(2^{-3},2^{-6},2^{-9},2^{-12})$ to get
\begin{equation}
(\mu^*_{4,10},
\mu^*_{5,10},
\mu^*_{6,10},
\mu^*_{7,10})<(5935,1,1,1).
\end{equation}
This information combines with 
Equation \ref{BOUND} and
Equation \ref{badbound} to give

\begin{equation}
|M_{J,10}|<19+5935+1+1+1<6000.
\end{equation}
For any real vector $V=(V_1,...,V_{343})$ define
\begin{equation}
\overline V=(|V_1|+6000,...,|V_{343}|+6000).
\end{equation}
This gives us
\begin{equation}
\sqrt 7 \epsilon F_{10}<\sqrt 7 \epsilon_0 \|\overline V_{10}\|<
1091<1234.
\end{equation}
This completes the proof of the Main Theorem
for $G^{\#}_{10}$.

\newpage

%% file: 7poly.tex
\section{Facts about Polynomials}

\subsection{Intervals and Interval Polynomials}
\label{rationalinterval}

We define a {\it rational interval\/} to be
an interval of the form $I=[L,R]$ where 
$L,R \in \Q$ and $L \leq R$.
For each operation $* \in \{+,-,\times\}$ we
define
\begin{equation}
I_1 * I_2=[\min(S),\max(S)],
\hskip 25 pt S=\{L_1*L_2,L_1*R_2,R_1*L_2,R_2*L_2\}.
\end{equation}
This definition is such that
$r_j \in I_j$ for $j=1,2$ then $r_1*r_2 \in I_1*I_2$.
Moreover, $I_1*I_2$ is the minimal
interval with this property.  The minimality
property implies that our laws are
both associative and distributive:
\begin{itemize}
\item $(I_1+I_2) \pm I_3=I_1+(I_2 \pm I_3)$.
\item $I_1 \times (I_2 \pm I_3)=(I_1 \times I_2) \pm (I_1 \times I_3)$.
\end{itemize}
We also can raise a rational interval to a 
nonnegative integer power:
\begin{equation}
I^k=I \times ... \times I \hskip 30 pt {\rm k\ times\/}.
\end{equation}

An {\it inteval polynomial\/} is an
expression of the form 
\begin{equation}
I_0+I_1t+ ... + I_n t^n.
\end{equation}
in which each of the coefficients are intervals
and $t$ is a variable meant to be taken in $[0,1]$.
Given the rules above, interval polynomials
may be added, subtracted or multiplied, in the
obvious way.  The associative and distributive
laws above give rise to similar results about
the arithmetic operations on interval polynomials.

We think of an ordinary polynomial as an interval
polynomial, just by taking the intervals to have
$0$ width.  We think of a constant as an interval
polynomial of degree $0$.  Thus, if we have some
expression which appears to involve constants,
ordinary polynomials, and interval polynomials,
we interpret everything in sight as an 
interval polynomial and then perform the
arithmetic operations needed to simplify
the expression.

Let $\cal P$ be the above interval polynomial.
We say that $\cal P$
{\it traps\/} the ordinary
polynomial 
\begin{equation}
C_0+C_1t+...+C_nt^n
\end{equation}
of the same degree
degree if $C_j \in I_j$ for all $j$.
We define the {\it min\/} of an interval polynomial
to be the polynomial whose coefficients are
the left endpoints of the intervals.  We
define the {\it max\/} similarly.   If $\cal P$
is an interval polynomial which traps an
ordinary polynomial, then 
\begin{itemize}
\item $\cal P$ traps ${\cal P\/}_{\rm min\/}$.
\item $\cal P$ traps ${\cal P\/}_{\rm max\/}$.
\item For all $t \in [0,1]$ we have
 ${\cal P\/}_{\rm min\/}(t) \leq P(t) \leq  {\cal P\/}_{\rm mat\/}(t)$.
\end{itemize}

Our arithmetic operations are such that if
the polynomial
${\cal P\/}_j$ traps the polynomial $P_j$ for
$j=1,2$, then
${\cal P\/}_1 * {\cal P\/}_2$ traps
$P_1*P_2$.  Here $* \in \{+,-,\times\}$.

\subsection{Rational Approximations of Power Combos}
\label{eig}

Suppose that $Y=(a_2,a_3,a_4,b_2,b_3,b_4)$ is a $6$-tuple
of rational numbers.   We are interested
in expressions of the form
\begin{equation}
C_Y(x)=
a_2  2^{-s/2} +
a_3  3^{-s/2} +
a_4  4^{-s/2} +
b_2 s 2^{-s/2} +
b_3 s 3^{-s/2} +
b_4 s 4^{-s/2}
\end{equation}
evaluated on the interval $[-2,16]$.
We call such expressions {\it power combos\/}.

For each even integer $2k=-2,...,16$.
we will construct rational polynomials
$A_{Y,2k,-}$ and $A_{Y,2k,+}$ and
$B_{Y,2k,-}$ and $B_{Y,2k,+}$ such that
\begin{equation}
\label{trap1}
A_{Y,2k,-}(t)  \leq C_Y(2k-t) \leq A_{Y,2k,+}(t), \hskip 30 pt
t \in [0,1].
\end{equation}
\begin{equation}
\label{trap2}
B_{Y,2k,-}(t)  \leq C_Y(2k+t) \leq B_{Y,2k,+}(t), \hskip 30 pt
t \in [0,1].
\end{equation}
We ignore the cases $(2,-)$ and $(16,+)$.

The basic idea is to use Taylor's Theorem with
Remainder:
\begin{equation}
\label{tayser}
m^{-s/2}=\sum_{j=0}^{11} \frac{(-1)^j \log(m)^j}{m^k 2^j j!} (s-2k)^j+
\frac{E_s}{12!} (s-2k)^{12}.
\end{equation}
Here $E_s$ is the
$12$th derivative of $m^{-s/2}$ evaluated at some point
in the interval.  Note that the only dependence on $k$ is
the term $m^k$ in the denominator, and this is a
rational number.

The difficulty with this approach is that
the coefficients of the above 
Taylor series are not rational. We get around
this trick by using interval polynomials.
We first pick specific intervals which
trap $\log(m)$ for $m=2,3,4$.  We choose
the intervals
$$
L_2=\bigg[\frac{25469}{36744},\frac{7050}{10171}\bigg], \hskip 10 pt
L_3=\bigg[\frac{5225}{4756},\frac{708784}{645163}\bigg], \hskip 10 pt
L_4=\bigg[\frac{25469}{18372},\frac{345197}{249007}\bigg].
$$
Each of these intervals has width about $10^{-10}$.
I found them using Mathematica's Rationalize function.
It is an easy exercise to check that
$\log(m) \in L_m$ for $m=2,3,4$.  It is also
an easy exercise to show that
\begin{equation}
\sup_{s \in [-2,16]} \max_{m=2,3,4}
\bigg|\frac{d^{12}}{ds^{12}} m^{-s/2} \bigg|<1.
\end{equation}
Indeed, the true answer is closer to $1/256$.
What we are saying is that we always have
$|E_s|<1$ in the series expansion from
Equation \ref{tayser}.  Fixing $k$ we introduce
the interval Taylor series

\begin{equation}
\label{tayser2}
A_{m}(t)=\sum_{j=0}^{11} \frac{(+1)^j (L_m)^j}{m^k 2^j j!}t^j+
\left[-\frac{1}{12!},\frac{1}{12!}\right] t^{12}.
\end{equation}

\begin{equation}
\label{tayser3}
B_m(t)=\sum_{j=0}^{11} \frac{(-1)^j (L_m)^j}{m^k 2^j j!}t^j+
\left[-\frac{1}{12!},\frac{1}{12!}\right] t^{12}.
\end{equation}

By construction $A_m$ traps the Taylor series
expansion from Equation \ref{tayser} when it
is evaluated at $t=s-2k$ and $t \in [0,1]$.
Likewise 
$B_m$ traps the Taylor series expansion 
from Equation \ref{tayser} when it is
evaluated at $t=2k-s$ and $t \in [0,1]$.
Define
$$
A_Y(t)=
a_2  A_2(t)+
a_3  A_3(t)+
a_4  A_4(t)+$$
\begin{equation}
b_2 (2k-t) A_2(t)+
b_3 (2k-t) A_3(t)+
b_4 (2k-t) A_4(t)
\end{equation}

$$
B_Y(t)=
a_2  B_2(t)+
a_3  B_3(t)+
a_4  B_4(t)+$$
\begin{equation}
b_2 (2k+t) B_2(t)+
b_3 (2k+t) B_3(t)+
b_4 (2k+t) B_4(t)
\end{equation}
By construction, $A_Y(t)$ traps
$C_Y(2k-t)$ when $t \in [0,1]$ and
$B_Y(t)$ traps $C_Y(2k+t)$ when
$t \in [0,1]$.

Finally, we define
$$
A_{Y,2k,-}=(A_Y)_{\rm min\/}, \hskip 15 pt
A_{Y,2k,-}=(A_Y)_{\rm max\/},$$
\begin{equation}
B_{Y,2k,-}=(B_Y)_{\rm min\/}, \hskip 15 pt
B_{Y,2k,-}=(B_Y)_{\rm max\/}.
\end{equation}
By construction these polynomials
satisfy Equations \ref{trap1} and
\ref{trap2} respectively for each $k=-1,...,8$.
These are our under and over approximations.
\newline
\newline
{\bf Remark:\/}
We implemented these polynomials in Java and
tested them extensively.

\subsection{Weak Positive Dominance}

Let 
\begin{equation}
\label{poly}
P(x)=a_0+a_1x+...+a_nx^n
\end{equation}
be a polynomial with real coefficients.  Here we describe
a method for showing that $P \geq 0$ on $[0,1]$,

Define
\begin{equation}
\label{sum1}
A_k=a_0+ \cdots + a_k.
\end{equation}
We call $P$ {\it weak positive dominant\/} (or
{\it WPD\/} for short) if
$A_k \geq 0$ for all $k$ and $A_n>0$.

\begin{lemma}
\label{PD1}
If $P$ is weak positive dominant, then
$P>0$ on $(0,1]$.
\end{lemma}

\startproof
The proof goes by induction on the degree of $P$.
The case $\deg(P)=0$ follows from the fact that $a_0=A_0>0$.
Let $x \in (0,1]$.
We have
$$P(x)=a_0+a_1x+x_2x^2+ \cdots + a_nx^n \geq $$
$$a_0x+a_1x+a_2x^2+ \cdots + a_nx^n=$$
$$x(A_1+a_2x+a_3x^2+ \cdots a_nx^{n-1})=xQ(x)> 0$$
Here $Q(x)$ is weak positive dominant and has degree $n-1$.
\endproof

Given an interval $I=[a,b] \subset \R$, let $A_I$ be one of
the two affine maps which carries $[0,1]$ to $I$. We
call the pair $(P,I)$ {\it weak positive dominant\/} if
$P \circ A_I$ is WPD. If $(P,I)$ is WPD
then $P \geq 0$ on $(a,b]$, by Lemma \ref{PD1}. 
For instance, if $P$ is WPD on $[0,1/2]$ and
$[1/2,1]$ then $P>0$ on $(0,1)$. 

\subsection{The Positive Dominance Algorithm}

Here I describe a method for certifying that
a polynomial (of several variables) is 
non-negative on a polytope.  I will
restrict to the case when
the polytope is the unit cube. I have used this
method extensively in other contexts.
See e.g. [{\bf S2\/}]. I don't know if
this method already exists in the literature.
It is something I devised myself.  

Given a multi-index 
$I=(i_1,...,i_k) \in (\N \cup \{0\})^k$ we
let 
\begin{equation}
x^I=x_1^{i_1}...x_k^{i_k}.
\end{equation}
Any polynomial $F \in \R[x_1,...,x_k]$
can be written succinctly as
\begin{equation}
F=\sum a_I X^I, \hskip 30 pt a_I \in \R.
\end{equation}
If $I'=(i_1',...,i_k')$ we write
$I' \leq I$ if $i'_j \leq i_j$ for
all $j=1,...,k$.
We call $F$ {\it positive dominant\/}
(PD)
if
\begin{equation}
\label{summa}
A_I:=\sum_{I' \leq I} a_{I'}> 0
\hskip 30 pt \forall I,
\end{equation}

\begin{lemma}
\label{PD2}
If $P$ is PD, then $P> 0$ on $[0,1]^k$.
\end{lemma}

\startproof
When $k=1$ the proof is the
same as in Lemma \ref{PD1},
once we observe that also $P(0)>0$.
Now we prove the general case.
Suppose the the coefficients of $P$ are
$\{a_I\}$.
We write
\begin{equation}
P=f_0+f_1x_k+...+f_mx_k^m,
\hskip 20 pt f_j \in \R[x_1,...,x_{k-1}].
\end{equation} 
Let $P_j=f_0+...+f_j$.
A typical coefficient in $P_j$ has
the form
\begin{equation}
\label{expandout}
b_J=\sum_{i=1}^j a_{Ji},
\end{equation}
where $J$ is a multi-index of length $k-1$
and $Ji$ is the multi-index of length $k$
obtained by appending $i$ to $J$.
From equation \ref{expandout} and
the definition of PD, the fact that
$P$ is PD implies that $P_j$ is
PD for all $j$.
\endproof

The positive dominance criterion is not that useful in
itself, but it feeds into a powerful divide-and-conquer
algorithm. We define the maps
$$
A_{j,1}(x_1,...,x_k)=(x_1,...,x_{j-1}\frac{x_j+0}{2},x_{i+1},...,x_k),
$$
\begin{equation}
\label{polysub}
A_{j,2}(x_1,...,x_k)=(x_1,...,x_{j-1}\frac{x_j+1}{2},x_{j+1},...,x_k),
\end{equation}

We define the $j$th {\it subdivision\/} of $P$ to be the
set
\begin{equation}
\label{SUB}
\{P_{j1},P_{j2}\}=
\{P \circ A_{j,1},P \circ A_{j,2}\}.
\end{equation}

\begin{lemma}
\label{half}
$P > 0$ on $[0,1]^k$ if and only if
$P_{j1}> 0$ and $P_{j2}> 0$ on
$[0,1]^k$.
\end{lemma}

\startproof
By symmetry, it suffices to take $j=1$.
Define
\begin{equation}
[0,1]^k_1=[0,1/2] \times [0,1]^{k-1}, \hskip 30 pt
[0,1]^k_2=[1/2,1] \times [0,1]^{k-1}.
\end{equation}
Note that
\begin{equation}
A_1([0,1]^k)=[0,1]^k_1,\hskip 30 pt
B_1 \circ A_1([0,1]^k)=[0,1]^k_2.
\end{equation}
Therefore,
$P > 0$ on $[0,1]^k_1$ 
if and only if $P_{j1} > 0$ on $[0,1]^k$. Likewise
$P > 0$ on $[0,1]^k_2$ if and only if
if $P_{j2} > 0$ on $[0,1]^k$.
\endproof

Say that a {\it marker\/} is a non-negative
integer vector in $\R^k$. 
Say that the {\it youngest entry\/} in the the marker
is the first minimum entry going from left to right. 
The {\it successor\/} of a marker is the marker obtained
by adding one to the youngest entry. For instance,
the successor of $(2,2,1,1,1)$ is $(2,2,2,1,1)$.
Let $\mu_+$ denote the successor of $\mu$.

We say that a {\it marked polynomial\/} is a pair
$(P,\mu)$, where $P$ is a polynomial and
$\mu$ is a marker.  Let $j$ be the position of the
youngest entry of $\mu$.  We define the
{\it subdivision\/} of $(P,\mu)$ to be the
pair
\begin{equation}
\{(P_{j1},\mu_+),(P_{j2},\mu_-)\}.
\end{equation}
Geometrically, we are cutting the domain in half
along the longest side, and using a particular rule
to break ties when they occur.
\newline
\newline
{\bf Divide-and-Conquer Algorithm:\/}
\begin{enumerate}
\item Start with a list LIST of marked polynomials.
Initially, LIST consists only of the marked polynomial
$(P,(0,...,0))$.
\item Let $(Q,\mu)$ be the last element of LIST.
We delete $(Q,\mu)$ from LIST and test whether
$Q$ is positive dominant.
\item Suppose $Q$ is positive dominant.
We go back to Step 2 if LIST is not empty.
Otherwise, we halt.
\item Suppose $Q$ is not positive dominant.
we append to LIST the two marked polynomials
in the subdivision of $(Q,\mu)$ and then go
to Step 2.
\end{enumerate}

If the algorithm halts, it constitutes a proof
that $P > 0$ on $[0,1]^k$.  Indeed, the
algorithm halts if and only if $P>0$ 
on $[0,1]^k$.
\newline
\newline
{\bf Parallel Version:\/}
Here is a variant of the algorithm.
Suppose we have a list $\{P_1,...,P_m\}$ of
polynomials and we want to show that at least
one of them is positive at each point of
$[0,1]^k$.  We do the following
\begin{enumerate}
\item Start with $m$ lists LIST$(j)$ for $j=1,...,m$
 of marked polynomials.
Initially, LIST$(j)$ consists only of the marked polynomial
$(P_j,(0,...,0))$.
\item Let $(Q_j,\mu)$ be the last element of LIST$(j)$.
We delete $(Q_j,\mu)$ from LIST$(j)$ and test whether
$Q_j$ is positive dominant.  We do this for $j=1,2,...$
until we get a success or else reach the last index.
\item Suppose {\it at least one\/}
$Q_j$ is positive dominant.
We go back to Step 2 if LIST$(j)$ is not empty.
(All lists have the same length.)
Otherwise, we halt.
\item Suppose none of $Q_1,...,Q_m$ is positive dominant.
For each $j$ we append to LIST$(j)$ the two marked polynomials
in the subdivision of $(Q_j,\mu)$ and then go
to Step 2.
\end{enumerate}
If this algorithm halts it constitutes a proof that
at least one $P_j$ is positive at each point
of $[0,1]^k$.

\subsection{Discussion}

For polynomials in $1$ variable, the
method of Sturm sequences counts
the roots of a polynomial in any
given interval. An early version of
this paper used Sturm sequences, but
I prefer the positive dominance criterion.
The calculations for the positive
dominance criterion are much simpler
and easier to implement.

There are generalizations of Sturm
sequences to higher dimensions, and
also other positivity criteria (such as
the Handelman decomposition) but
I bet they don't work as well as the
positive dominance algorithm.
Also, I don't see how to do 
the parallel positive dominance
algorithm with these other methods.

The positive dominance algorithm works
so well that one can ask why I didn't
simply use it to prove the Main Theorem
straight away.  After all, the Main
Theorem does reduce to a positivity
theorem about a finite set of polynomials.
I tried this.  However, the polynomials
seem to involve an astronomical number
of terms.  It is not a feasible calculation.

\newpage

%% file: 8tumanov.tex
\section{Proof of Lemma \ref{Tumanov}}

\subsection{Some General Considerations}
\label{general}

First we discuss the general principle behind
Lemma \ref{Tumanov}.  Suppose that
the TBP is a minimizer with respect
to $\Gamma_1$ and $\Gamma_2$, and a unique minimizr
with respect to
$\Gamma_3,...,\Gamma_m$.
Suppose also that $R$ is some other energy
function and we want to show that the TBP
is the unique minimizer with respect
to $R$.  This is true if we can find
a combination
\begin{equation}
\Gamma=a_0+a_1\Gamma_1+...+a_m\Gamma_m,
\hskip 30 pt a_1,...,a_k \geq 0.
\end{equation}
Here $a_0$ could be negative. This doesn't bother us.
such that
\begin{itemize}
\item The constants $a_3,...,a_k$ do not identically vanish.
\item $\Gamma(x) \leq R(x)$ for all $x \in (0,2]$.
\item $\Gamma(x)=R(x)$ for $x=\sqrt 2,\sqrt 3,\sqrt 4$.
\end{itemize}
The three special values are the distances between
pairs of points of the TBP.  We find it nice to write
$\sqrt 4$ instead of $2$.

Let $X_0$ be the TBP and let $X$ be any other configuration
of $5$ distinct points on the sphere.  Letting
$\Gamma_0$ be the function which is identically $1$,
we have
$$
{\cal E\/}_R(X) \geq
{\cal E\/}_{\Gamma}(X)=
\sum_{i=0}^m a_i {\cal E\/}_{\Gamma_k}(X)
>\sum_{i=0}^m a_i {\cal E\/}_{\Gamma_k}(X_0)=
{\cal E\/}_{\Gamma}(X_0)={\cal E\/}_R(X_0).
$$

Here is how we find such positive combinations.
We set $m=4$, so that we are looking for $5$
coefficients $a_0,...,a_4$.  We impose the
$5$ conditions
\begin{itemize}
\item $\Gamma(x)=R(x)$ for $x=\sqrt 2,\sqrt 3,\sqrt 4$.
\item $\Gamma'(x)=R'(x)$ for $x=\sqrt 2,\sqrt 3$.
\end{itemize}
Here $R'=dR/dx$ and $\Gamma'=d\Gamma/dx$.
These $5$ conditions give us $5$ linear
equations in $5$ unknowns.
In the cases described below, the
associated matrix is invertible and there is a unique
solution.  In our situation it will be obvious that
the constants $a_1,a_2,a_3,a_4$ cannot identically
vanish.  That $\Gamma \leq R$ is far from
obvious, but the techniques from the previous chapter
will estblish this in each case of interest.

\subsection{Finding the Coefficients}

Recall that $R_s(r)=r^{-s}$ when $s>0$ and
$R_s(r)=-r^{-s}$ when $s<0$.  
We break Lemma \ref{Tumanov} into $3$ cases.

\begin{enumerate}
\item When $s \in (-2,0)$ we
use $G_1,G_2,G_3,G_5$.
\item When $s \in (0,6]$ we
use $G_1,G_2,G_4,G_6$.
\item When $s \in [6,13]$ we
use $G_1,G_2,G_5, G^{\#}_{10}$.
\end{enumerate}

We will also keep track of the expression
\begin{equation}
\delta=2\Gamma'(2)-2R'(2).
\end{equation}

In the first case we have
$$
\left[\matrix{
a_0 \cr a_1 \cr a_2 \cr a_3 \cr a_4\cr \delta}\right]=
\frac{1}{144} \left[\matrix{
0 & 0 & -144 & 0 & 0 & 0 \cr 
-312 & -96 & 408 & 24 & 80 & 0 \cr 
684 & -288 & -396 & -54 & -144 & 0 \cr 
-402 & 264 & 138 & 33 & 68 & 0 \cr 
30 & -24 & -6 & -3 & -4 & 0 \cr 
2496 & 768 & -3264 & -192 & -640 & -144}\right]
 \left[\matrix{
2^{-s/2}\cr 3^{-s/2}\cr 4^{-s/2} \cr
s2^{-s/2}\cr s3^{-s/2}, \cr s4^{-s/2}}\right].
$$

In the second case we have
$$
\left[\matrix{
a_0 \cr a_1 \cr a_2 \cr a_3 \cr a_4\cr \delta}\right]=
\frac{1}{792}\left[\matrix{
0 & 0 & 792 & 0 & 0 & 0 \cr
792 & 1152 & -1944 & -54 & -288 & 0 \cr 
-1254 & -96 & 1350 & 87 & 376 & 0 \cr
528 & -312 & -216 & -39 & -98 & 0 \cr
-66 & 48 & 18 & 6 & 10 & 0 \cr
-6336 & -9216 & 15552 & 432 & 2304 & 792}\right]
 \left[\matrix{
2^{-s/2}\cr 3^{-s/2}\cr 4^{-s/2} \cr
s2^{-s/2}\cr s3^{-s/2}, \cr s4^{-s/2}}\right].
$$

In the third case we have
{\footnotesize
$$
\left[\matrix{
a_0 \cr a_1 \cr a_2 \cr a_3 \cr \widehat a_4\cr \delta}\right]=
\frac{1}{268536}\left[\matrix{
0& 0& 268536& 0& 0& 0 \cr 
88440& 503040& -591480& -4254& -65728& 0 \cr 
-77586& -249648& 327234& 2361& 65896& 0\cr 
41808& -19440& -22368& -2430& -9076& 0\cr 
-402& 264& 138& 33& 68& 0\cr 
-707520& -4024320& 4731840& 34032& 525824& 268536}\right] 
\left[\matrix{
2^{-s/2}\cr 3^{-s/2}\cr 4^{-s/2} \cr
s2^{-s/2}\cr s3^{-s/2}, \cr s4^{-s/2}}\right].
$$
}

Thus the coefficients are precisely the power combos
considered in \S \ref{eig}.  These power combos
are functions of the variable $s$.

\subsection{Positivity Proof}

Now we explain how we prove that
$a_1,a_2,a_3,a_4,\delta>0$ on the
relevant intervals.
We will do the three cases one at a time.
We consider $a_1$ on $(-2,0)$ in detail.

\begin{itemize}
\item We set
$Y=(-312,-96,408,24,80,0)$, the
row of the relevant matrix corresponding to $a_1$.
By construction,
$a_1(s)=C_Y(s)$.
\item We verify that the two
under-approximations
$A_{Y,-2,+}$ and
$A_{Y,0,-}$ are WPD on 
$[0,1/2]$ and on $[1/2,1]$.
Hence these functions are
positive on $(0,1]$. 
\item Since 
$A_{Y,-2,+}(t) \leq a_1(t-2)$ for $t \in [0,1]$ we
see that $a_1>0$ on $(-2,-1]$.

\item Since 
$A_{Y,0,-}(t) \leq a_1(-t)$ for $t \in [0,1]$ we
see that $a_1>0$ on $[1,0)$.
\end{itemize}

The same argument works for $a_2,a_3,a_4,\delta$ on
$[-2,0]$.   In each case, the relevant
under-approximation is either WPD
on $[0,1]$ or WPD on $[0,1/2]$ and
$[1/2,1]$. 
We conclude that $a_1,a_2,a_3,a_4,\delta>0$ on
$(-2,0)$.
The statement 
$\delta>0$ means that $\Gamma'(2)>R'(2)$ for all
$s \in (-2,0)$.  
For later use we use the same method to check
that $$\Gamma(0)=c_0+4c_1+16c_2+256c_3+1024c_5<0$$ for
all $s \in (-2,0)$.

We do the same thing on the interval $[0,6]$, except
that we use the intervals
$[0,1]$, $[1,2]$, $[2,3]$, $[3,4]$, $[4,5]$ and $[5,6]$
and the corresponding under-approximations.
In this case, every polynomial in sight -- all
$30=5 \times 6$ of them -- is WPD
on $[0,1]$.  We conclude that
that $a_1,a_2,a_3,a_4,\delta>0$ on $(0,6]$.
To be sure, we check the interval endpoints
$s=1,2,3,4,5,6$ by hand.

Finally, we do the same thing on the interval $[6,13]$,
using the intervals $[6,7],...,[12,13]$.
Again, every polynomial in sight is WPD on $[0,1]$,
and in fact PD on $[0,1]$.   Since we just check
the WPD condition, we also check that our
functions are positive at the integer values in
$[6,13]$ by hand.
We conclude 
that $a_1,a_2,a_3,a_4,\delta>0$ on $[6,13]$.
\newline
\newline
{\bf Remark:\/}
In the third case, we checked additionally
that $a_1,a_2,a_3,a_4,\delta$ are positive
on $[3,13+1/16]$. Thus, our result is really
true for all power law exponents up to
$13+1/16$.  The reader can play with our
graphical user interface, see plots of all
these functions, and run positivity tests.

\subsection{Under Approximation: Case 1}
\label{ua}

Here we show that
$\Gamma_s(r) \leq R_s(r)$ for all $r \in (0,2]$
and all $s \in (-2,0)$.
We suppress the dependence on $s$ as much
as we can.  In particular, we set $R=R_s$, etc.
Here $R<0$ so we want $\Gamma/R>1$.
Define
\begin{equation}
\label{under}
H(r)=\frac{\Gamma}{R}-1=-r^s \Gamma-1.
\end{equation}
We just have to show that $H \geq 0$ on $(0,2)$.
Let $H'=dH/dr$.

\begin{lemma}
\label{Qdefined}
$H'$ has $4$ simple roots in $(0,2)$.
\end{lemma}

\startproof
We count roots with multiplicity.
We have
\begin{equation}
\label{deriv}
H'(r)= - r^{s-1}(s\Gamma(r)+r\Gamma'(r)).
\end{equation}
Comgining Equation \ref{deriv} with the general 
equation
\begin{equation}
rG_k'(r)=2kG_k(r)-8kG_{k-1}(r),
\end{equation}
we see that 
the positive roots of $H'(r)$ are the same as
the positive roots of
\begin{equation}
- r^{s-1}H'(r)=(10+s)c_4 G_5(r)-40 c_4 G_4(r)+\sum_{k=1}^3 b_k G_k(r)+b_0.
\end{equation}
Here $b_0,...,b_3$ are coefficients we don't care about.
Making the substitution $t=4-r^2$ we see that the roots of
$H'$ in $(0,2)$ are in bijection with the roots in $(0,4)$ of
\begin{equation}
\label{monic}
\psi(t)=t^5-\frac{40}{10+s} t^4 + b_3t^3+b_2t^2+b_1t+b_0.
\end{equation}
Moreover, the change of coordinates is a diffeomorphism
from $(0,4)$ to $(0,2)$ and so it carries simple
roots to simple roots.

The polynomial $\psi$ has $5$ roots counting 
multiplicities.  Let's find $4$ of these roots
first. Since $H(\sqrt 2)=H(\sqrt 3)=H'(\sqrt 2)=H'(\sqrt 3)=0$,
we see that $H'$ has at least $4$ roots in $(0,2)$.
Besides the roots at $\sqrt 2$ and $\sqrt 3$, there
is a root in $(\sqrt 2,\sqrt 3)$ and a root in $(\sqrt 3,2)$.
This means that $\psi$ has $4$ corresponding
roots in $(0,4)$.  We claim that $\psi$ has an
even number of roots in $(0,2)$, counting multiplicity.
Once we know this, we can say that
the $4$ roots we have found are simple.
But then the corresponding $4$ roots of
$H'$ in $(0,2)$ are simple and there are no
others.

Now for the parity argument.
Since $R<0$ on $(0,2)$ and
$R(0)=0$ and $\Gamma(0)<0$ we see that
$H(r) \to \infty$ as $r \to 0$.  
Hence $\psi(t) \to \infty$ as $t \to 4$.
Hence, there are arbitrarily small values
of $\delta>0$ such that $\psi(4-\delta)>0$.

Since
$\Gamma'(2)>R'(2)$ and $R(2)<0$ we see that
$H'(2)<0$.  Since our
change of coordinates is an
orientation reversing diffeomorphism
we have $\psi(0)>0$ by the chain rule.

Now we know that there are arbitrarily
small values $\delta>0$ so that
$\psi(\delta)>0$ and $\psi(4-\delta)>0$.
But then the number of roots of
$\psi$ in $(\delta,4-\delta)$ is even.
Since $\delta$ is arbitrary, the number
roots of $\psi$ in $(0,4)$ is even.
\endproof

\begin{lemma}
$H''(\sqrt 2)>0$ and
$H''(\sqrt 3)>0$ for all 
$s \in (-2,0)$.
\end{lemma}

\startproof
How we mention the explicit dependence on 
$s$ and remember that we are taking about $H_s$.
We check directly that $H_{-1}''(\sqrt 2)>0$.
It cannot happen that $H_s''(\sqrt 2)=0$
for other $s \in (-2,0)$ 
because then $H'_s$ does not have only
simple roots in $(0,2)$.  Hence
$H_s''(\sqrt 2)>0$ for all $s \in (-2,0)$.
The same argument shows that
$H_s''(\sqrt 3)>0$ for all $s \in (-2,0)$.
\endproof

Now we set $H=H_s$ again.

\begin{lemma}
For all sufficiently small $\delta>0$ the
quantities
$$H(0+\delta), \hskip 20 pt
H(\sqrt 2 \pm \delta), \hskip 20 pt
H(\sqrt 3 \pm \delta), \hskip 20 pt
H(2-\delta)$$
are positive.
\end{lemma}

\startproof
We have already seen that $H(\delta) \to +\infty$
as $\delta \to 0$.  Likewise, we have seen that
$H'(2)<0$ and $H(2)=0$.  So $H(2-\delta)>0$ for
all sufficiently small $\delta$.
Finally, the case of $\sqrt 2$ and $\sqrt 3$
follows from the previous lemma and the second
derivative test.
\endproof

We already know that $H'$ has exactly
$4$ simple roots in $(0,2)$.  
In particular, the interval
$(0,\sqrt 2)$ has no roots of $H'$
and the intervals $(\sqrt 2,\sqrt 3)$
and $(\sqrt 3,2)$ have $1$ root each.
Finally, we know that $H>0$ sufficiently
near the endpoints of all these intervals.
If $H(x)<0$ for some $x \in (0,2)$, then $x$ must
be in one of the $3$ intervals just mentioned,
and this interval contains at least $2$
roots of $H'$.  This is a contradiction.

\subsection{Under Approximation: Case 2}

This time we have $s \in (0,6]$.
We have $R>0$ so we set
\begin{equation}
H=\frac{\Gamma}{R}-1.
\end{equation}
It suffices to prove that
$H \geq 0$ on $(0,2]$.

Everything in Case 1 works here, word for word,
provided that Lemma \ref{Qdefined} holds
for $H$, when $s \in (0,6]$.  The proof
is exactly the same, except for a global
sign, and the fact that this time we have

\begin{equation}
\label{monic2}
\psi(t)=t^6-\frac{48}{12+s} t^4 + b_4t^4 b_3t^3+b_2t^2+b_1t+b_0.
\end{equation}

We just have to show that $\psi$ has $4$ simple
roots in $(0,4)$.  Note that the
sum of the $6$ roots of $\psi$ is
$48/(12+s)<4$.  This works because $s>0$ here.
The $4$ roots of $\psi$ we already know
about are $1$ and $2$ and some number
in $(0,1)$ and some number in $(1,2)$.
The sum of these roots exceeds $4$ and
so the remaining two roots cannot also
be positive.  Hence $\psi$ has at most
$5$ roots in $(0,4)$.

The parity argument works the same way.
This time $H(\delta) \to \infty$ as $\delta \to 0$
because $R \to \infty$ and $\Gamma$ is bounded.
The parity argument shows that
$\psi$ has an even number of roots in $(0,4)$.
Hence, $\psi$ has exactly $4$ such roots
and they are all simple.  Hence $H'$ has
exactly $4$ simple roots in $(0,2)$.
This completes the proof in Case 2.

\subsection{Under Approximation: Case 3}
\label{case3}

This time we have $s \in [6,16]$ and
everything is as in Case 2.
All we have to do is show that the polynomial
$\psi$, as Equation \ref{monic2}, has
exactly $4$ simple roots in $(0,2)$.
This time we have
\begin{equation}
\label{monic3}
\psi(t)=t^{10}-\frac{80}{20+s} t^9 + b_8t^8+...+b_0.
\end{equation}
Again these coefficients depend on $s$.
\newline
\newline
{\bf Remark:\/}
$\psi$ only has $7$ nonzero terms and hence
can only have $6$ positive real roots.
The number of positive real roots is bounded
above by Descartes' rule of signs.
Unfortunately, $\psi$ turns
out to be alternating and so Descartes' rule
of signs does not eliminate the case of $6$ roots.
This approach seems useless.
The sum of the roots of $\psi$
is less than $4$, so it might seem as if
we could proceed as in Case 2.  Unfortunately, there are
$10$ such roots and this approach also seems useless.
We will take another approach to proving what we want.
\newline

\begin{lemma}
When $s=6$ the polynomial $\psi$ has $4$ simple
roots in $(0,4)$.
\end{lemma}

\startproof
We compute explicitly that
$$
\psi(t)=t^{10}-\frac{40}{13}t^9+ 
\frac{830304}{5785}t^5 -
\frac{415152}{1157} t^4 +\frac{789255}{1157} t^2
-\frac{3264104}{5785} t + \frac{115060}{1157}.
$$
This polynomial only has $4$ real roots -- the
ones we know about.  The remaining roots are
all at least $1/2$ away from the interval
$(0,4)$ and so even a very crude analysis
would show that these roots do not lie
in $(0,4)$.  We omit the details.

\begin{lemma}
Suppose, for all $s \in [6,13]$ that
$\psi$ only has simple roots in
$(0,4)$.  Then in all cases
$\psi$ has exactly $4$ such roots.
\end{lemma}

\startproof
Let $N_s$ denote the number of simple roots
of $\psi$ at the parameter $s$.  The same
argument as in Cases 1 and 2 shows that
$N_s$ is always even.  Suppose $s$ is
not constant.  Consider the infimal
$u \in (6,13]$ such that $N_u>4$.
The roots of $\psi$ vary continuously with
$s$.  How could more roots move into
$(0,4)$ as $s \to u$?

One possibility is that such a
root $r_s$ approaches from the upper half
plane or from the lower half plane.
That is, $r_s$ is not real for $s<u$.
Since $\psi$ is a real polynomial,
the conjugate $overline r_s$ is also
a root.  The two roots
$r_s$ and $\overline r_s$ are approaching
$(0,4)$ from either side. But then the the limit
$$\lim_{s \to u} r_s$$ is a double root of
$\psi$ in $(0,4)$.  This is a contradiction.

The only other possibility is that the
roots approach along the real line.
Hence, there must be some $s<u$ such that
both $0$ and $4$ are roots of $\psi$.
But the same parity argument as in Case 1
shows $\psi(0)>0$ for all $s \in [6,13]$. 
\endproof

To finish the proof we just have to show that
$\psi$ only has simple roots in $(0,4)$ for
all $s \in [6,13]$.  We bring the
dependence on $s$ back into our notation and
write $\psi_s$.  It suffices to show that
that $\psi_s$ and $\psi'_s=d\psi_s/dr$ do not simultaneously
vanish on the rectangular domain
$(s,r) \in [6,13] \times [0,4]$.
This is a job for our method of positive
dominance.

We will explain in detail what do on the
smaller domain $$(s,r)=[6,7] \times [0,4].$$
The proof works the same for the remining
$1 \times 4$ rectangles.
The coefficients of $\psi_s$ and $\psi'_s$ are
power combos in the sense of
Equation \ref{eig}.  

We have rational vectors 
$Y_0,...,Y_9$ such that
\begin{equation}
\psi_s(r)=
\sum_{j=0}^9 C_j r^j, \hskip 20 pt
C_j=C_{Y_j}.
\end{equation}

We have the under- and over-approximations:
\begin{equation}
A_j=A_{Y_j,6,+}
\hskip 30 pt
B_j=B_{Y_j,6,+}
\end{equation}

We then define $2$-variable under- and over-approximations:
\begin{equation}
\underline \psi(t,u)= \sum_{i=0}^9 A_i(t) (4u)^i,
\hskip 30 pt
\overline \psi(t,u)= \sum_{i=0}^9 B_i(t) (4u)^i.
\end{equation}

We use $4u$ in these sums because we want our
domains to be the unit square.  We have
\begin{equation}
\label{overunder}
\underline \psi(t,u) \leq 
\psi_{6+t}(4u) \leq
\overline \psi(t,u), 
\hskip 30 pt
\forall (t,u) \in [0,1]^2.
\end{equation}

Now we do the same thing for $\psi'_s$.
We have rational vectors
$Y_0',...,Y_8' \in \Q^8$ which
work for $\psi'$ in place of $\psi$, and
this gives under- and over-approximations
$\underline \psi'$ and $\overline \psi'$ which
satisfy the same kind of equation as
Equation \ref{overunder}.

We run the parallel positive dominance algorithm
on the set of functions
$\{\underline \phi,-\overline \phi,
\underline \phi',-\overline \phi'\}$
and the algorithm halts.  This constitutes
a proof that at least one of these
functions is positive at each point.
But then at least one of $\phi_s(r)$
or $\phi'_s(r)$ is nonzero for each 
$s \in [6,13]$ and each $r \in [0,4]$.
Hence $\psi_s$ only has simple roots
in $(0,4)$.  This completes our proof
in Case 3.

Our proof of Lemma \ref{Tumanov} is done.

\newpage

%% file: 9compute.tex
\section{Computational Details}
\label{computer}

\subsection{Getting the Program}
\label{download}

Our computer program is written in Java. At least in 2016,
one can get it from my Brown University website:
\newline
\newline
{\bf http://www.math.brown.edu/$\sim$res/Java/Riesz.tar\/}
\newline
\newline
The directory has $3$ relevant subdirectories:
\begin{itemize}
\item {\bf Approximations:\/}
This has the computer code for
with Lemma \ref{Tumanov}.
\item {\bf Hessian:\/} This has the
Mathematica code used in the the
local analysis of the Hessian.
\item {\bf Riesz:\/}
This has the main program, which runs the
divide and conquer algorithm from \S 5.
\end{itemize}
This main directory has a README file which
contains more information about these
directories.  Each subdirectory has a
README file as well, which gives information
about running the program. The
{\bf Approximations\/} and {\bf Riesz\/} directories
each contain the code for java programs.
The programs each have a documentation window
which explains how to operate the program.
Each of the programs also has a debugging
mode, where the main operations are checked.

\subsection{Debugging}

One serious concern about any computer-assisted
proof is that some of the main steps of
the proof reside in computer programs which
are not printed, so to speak, along with
the paper.  It is difficult for one to
directly inspect the code without a serious
time investment, and indeed the interested
reader would do much better simply to reproduce
the code and see that it yields the same results.

The worst thing that could happen is
if the code had a serious bug which caused it
to suggest results which are not actually true.
Let me explain the extent to which I have debugged
the code.  Each of the java programs has a
debugging mode, in which the user can test that
various aspects of the program are running correctly. 
While the debugger does not check every method, it
does check that the main ones behave exactly
as expected.

Here is what the user can check in the
{\bf Approximations\/} program:
\begin{itemize}
\item Our under-approximation (respectively over-approximation)
of $\log(m)$ is less than (respectively greater than)
the numerically computed value of $\log(m)$. Here
$m=2,3,4$.
\item The series under-approximation
(respectively over-approximation) to a random
power combo 
evaluated at a random point in a random
unit integer interval is less than (respectively 
greater than) and
very close to the
power combo when it is computed numerically using
Java's power function.
\item The series under-approximations to
the functions $\psi$ and $\psi'$ from \S \ref{case3}
behave graphically as they should.  For instance,
when $\underline \psi$ is plotted alongside the graph of
the function $H=(1-\Gamma/F)$, the zeros of
$\underline \psi$ visually match the locations
of the extrema of $H$.  Likewise the zeros of
$\underline \psi'$ visually match the extrema
of $\underline \psi$.
\item The polynomial subdivision from
Equation \ref{polysub} checks out correctly
on random inputs -- both in the $1$-variable
case and in the $2$-variable case.
\end{itemize}

Here are the things one can check for the
main program, the one in the {\bf Riesz\/}
directory.

\begin{itemize}
\item You can check on random inputs that
the interval arithmetic operations are
working properly.
\item You can check on random inputs that
the vector operations - dot product, addition, etc. - are
working properly.
\item You can check for random dyadic squares that
the floating point and interval arithmetic
measurements match in the appropriate sense.
\item You can select a block of your choice
and compare the estimate from the Energy
Theorem with the minimum energy taken over
a million random configurations in the block.
\item You can open up an auxiliary window
and see the grading step of the algorithm
performed and displayed for a block of
your choosing.
\end{itemize}

I have not included a debugger for the code
in the {\bf Hessian\/} directory because
this Mathematica code is quite straightforward.  
The code performs the calculations
described in detail in \S \ref{local}.

In any case, I view Lemmas \ref{energy1} and
\ref{energy2} as the main mathematical
contributions to the paper, and these
are covered by the code in the {\bf Riesz\/} directory.

\subsection{Integer Calculations}

Let me discuss the implementation of the
divide and conquer algorithm.
We manipulate blocks and dyadic squares
using {\bf longs\/}.  These are $64$ bit
integers.  Given a dyadic square $Q$
with center $(x,y)$ and side length $2^{-k}$,
we store the triple 
\begin{equation}
(Sx,Sy,k).
\end{equation}
Here $S=2^{25}$ when we do the calculations
for $G_3,G_4,G_5,G_6$ and
$S=2^{30}$ when we do the calculation for
$\widehat G_{10}$.  The reader can modify
the program so that it uses any power of 
$2$ up to $2^{40}$.  Similarly,
we store a dyadic segment with
center $x$ and side length
$2^{-k}$ as $(Sx,k)$.

The subdivision is then obtained by manipulating
these triples.  For instance, the top right square
in the subdivision of $(Sx,Sy,k)$ is
$$\big(Sx-2^{-k+1}S,Sy-2^{-k+1}S,k+1\big).$$
The scale $2^N$ allows for $N$ such subdivisions
before we lose the property that the squares
are represented by integer triples.
The biggest dyadic square is stored as
$(0,0,-2)$, and each subdivision increases
the value of $k$ by $1$.  We terminate the
algorithm if we ever arrive at a dyadic
square whose center is not an even pair of integers.
We never reach this eventuality when we
run the program on the functions from the
Main Theorem, but it does occur if we try
functions like $G_7$.

We make exact comparisons for Steps 1-4 in the
grading part of the algorithm described
in \S \ref{GRADE}. For instance, the point
of scaling our square centers by $S$ is that
the inequalities which go into the calculations
in \S \ref{checkTBP} are all integer inequalities.
We are simply clearing denominators.
It is only Step 5
which requires floating point (or interval)
calculations.

\subsection{Interval Arithmetic}

We implement interval arithmetic the same way that
we did in [{\bf S1\/}].  Here we repeat some
of the discussion, but abbreviate it.  We make
some changes to the way we do things in [{\bf S1\/}],
and also things are simpler here because we never
use the square root function.

Java represents real numbers by {\bf doubles\/}, essentially
according to the scheme
discussed in [{\bf I\/}, \S 3.2.2].  A double is a
$64$ bit string where $11$ of the bits control the
exponent, $52$ of the bits control the binary expansion,
and one bit controls the sign.
The non-negative doubles have a lexicographic ordering, and this ordering
coincides with the usual ordering of the real numbers they
represent.   The lexicographic ordering for the non-positive doubles
is the reverse of the usual ordering of the real numbers they
represent.  To {\it increment\/} $x_+$ of
 a positive double $x$ is the very next double
in the ordering.  This amounts to treating the last $63$ bits of the
string as an integer (written in binary) and adding $1$ to it.
With this interpretation, we have $x_+=x+1$.
We also have the decrement $x_-=x-1$.
Similar operations are defined on the non-positive doubles.
These operations are not defined on the largest and smallest
doubles, but our program never encounters (or comes anywhere near)
these.

Let $\D$ be the set of all doubles.
Let 
\begin{equation}
\R_0=\{x \in \R|\ |x| \leq 2^{50}\}
\end{equation}
Our choice of $2^{50}$ is an arbitrary but convenient cutoff.
Let $\D_0$ denote the set of doubles representing reals in $\R_0$.

According to the discussion in 
[{\bf I\/}, 3.2.2, 4.1, 5.6], there is a map
$\R_0 \to \D_0$ which maps each $x \in \R_0$ to some
$[x] \in \D_0$ which is closest to $x$.   In case there
are several equally close choices, the computer chooses one
according to the method
in [{\bf I\/}, \S 4.1].    This ``nearest point projection''
exists on a subset of $\R$ that is much larger
than $\R_0$, but we only need to consider
$\R_0$.  We also have the inclusion $r: \D_0 \to \R_0$, which
maps a double to the real that it represents.   

Our calculations just use the arithmetic
operations (plus, minus, times, divide)
These operations act on $\R_0$ in the usual way.
Operations with the same name act on $\D_0$.
Regarding these operations, [{\bf I\/}, \S 5] states that
{\it each of the operations shall be performed as if it first produced an
intermediate result correct to infinite precision and with unbounded range, and then 
coerced this intermediate result to fit into the destination's format\/}.  Thus,
for doubles $x,y \in \D_0$ such that $x*y \in \D_0$ as well, we have
\begin{equation}
\label{rule}
x*y=[r(x)*r(y)]; \hskip 20 pt
* \in \{+,-,\times,\div\}.
\end{equation}
The operations on the left hand side represent operations on doubles
and the operations on the right hand side represent operations on
reals.

It might happen that $x,y \in \D_0$ but $x*y$ is not.
To avoid this problem, we make the following
checks before performing any arithmetic operation.
\begin{itemize}
\item For addition and subtraction,
$\max(|x|,|y|)\leq2^{40}$.
\item For multiplication,
either $|x|\leq2^{40}$ and $|y|\leq2^{10}$ or
$|x|\leq2^{10}$ and $|y|\leq2^{40}$.
\item For division,
$|x|\leq2^{40}$ and $|y|\leq2^{10}$ and
$|y| \geq 2^{-10}$.
\end{itemize}
We set the calculation to abort if any of these
conditions fails.
\newline

For us, an
{\it interval\/} is a pair $I=(x,y)$ of doubles with $x \leq y$ and
$x,y \in \D_0$.   Say
that $I$ {\it bounds\/} $z \in \R_0$ if $x \leq [z] \leq y$.  This is true
if and only if $x \leq z \leq y$.  Define
\begin{equation}
[x,y]_o=[x_-,y_+].
\end{equation}
This operation is well defined for doubles in $\D_0$.
We are essentially {\it rounding out\/} the endpoints of the interval.
Let $I_0$ and $I_1$ denote the left and right endpoints of $I$.
Letting $I$ and $J$ be intervals, we define
\begin{equation}
\label{operate}
I*J = (\min_{ij} I_i *I_j,\max_{ij} I_i*I_j)_o.
\end{equation}

That is, we perform the operations on all the endpoints, order
the results, and then round outward. 
Given Equation \ref{rule}, we the interval
$I*J$ bounded $x*y$ provided that $I$ bounds $x$ and $J$ bounds $y$.
Except for the rounding out, this is the same
as what we discussed for the rational
interals in \S \ref{rationalinterval}.

We also define an interval version of a vector in $\R^3$.
Such a vector consists of $3$ intervals.  The only
operations we perform on such objects are addition, subtraction,
scaling, and taking the dot product.  These operations are
all built out of the arithmetic operations.

All of our calculations come down to proving inequalities 
of the form $x<y$.  We imagine that $x$ and $y$ are the
outputs of some finite sequence of arithmetic operations
and along the way we have intervals
$I_x$ and $I_y$ which respectively bound $x$ and $y$.
If we know that the right endpoint of $I_x$ is less
than the left endpoint of $I_y$, then this constitutes
a proof that $x<y$.  The point is that the whole
interval $I_x$ lies to the left of $I_y$ on the number line.

\newpage

%% file: refs.tex
\section{References}

[{\bf A\/}] A. N. Andreev,
{\it An extremal property of the icosahedron\/}
East J Approx {\bf 2\/} (1996) no. 4 pp. 459-462
\newline
\newline
[{\bf BBCGKS\/}] Brandon Ballinger, Grigoriy Blekherman, Henry Cohn, Noah Giansiracusa, Elizabeth Kelly, Achill Schurmann, \newline
{\it Experimental Study of Energy-Minimizing Point Configurations on Spheres\/}, 
arXiv: math/0611451v3, 7 Oct 2008
\newline
\newline
[{\bf BDHSS\/}] P. G. Boyvalenkov, P. D. Dragnev, D. P. Hardin, E. B. Saff, M. M. Stoyanova,
{\it Universal Lower Bounds and Potential Energy of Spherical Codes\/}, 
Constructive Approximation 2016 (to appear)
\newline
\newline
[{\bf C\/}] Harvey Cohn, {\it Stability Configurations of Electrons on a Sphere\/},
Mathematical Tables and Other Aids to Computation, Vol 10, No 55,
July 1956, pp 117-120.
\newline
\newline
[{\bf CK\/}] Henry Cohn and Abhinav Kumar, {\it Universally 
Optimal Distributions of Points on Spheres\/}, J.A.M.S. {\bf 20\/} (2007) 99-147
\newline
\newline
[{\bf CCD\/}] online website: \newline
http://www-wales.ch.cam.ac.uk/$\sim$ wales/CCD/Thomson/table.html
\newline
\newline
[{\bf DLT\/}] P. D. Dragnev, D. A. Legg, and D. W. Townsend,
{\it Discrete Logarithmic Energy on the Sphere\/}, Pacific Journal of Mathematics,
Volume 207, Number 2 (2002) pp 345--357
\newline
\newline
[{\bf HS\/}], Xiaorong Hou and Junwei Shao,
{\it Spherical Distribution of 5 Points with Maximal Distance Sum\/}, 
arXiv:0906.0937v1 [cs.DM] 4 Jun 2009
\newline
\newline
[{\bf I\/}] IEEE Standard for Binary Floating-Point Arithmetic
(IEEE Std 754-1985)
Institute of Electrical and Electronics Engineers, July 26, 1985
\newline
\newline
[{\bf RSZ\/}] E. A. Rakhmanoff, E. B. Saff, and Y. M. Zhou,
{\it Electrons on the Sphere\/}, \newline
  Computational Methods and Function Theory,
R. M. Ali, St. Ruscheweyh, and E. B. Saff, Eds. (1995) pp 111-127
\newline
\newline
[{\bf S1\/}] R. E. Schwartz, {\it The $5$ Electron Case of Thomson's Problem\/},
Journal of Experimental Math, 2013.
\newline
\newline
[{\bf S2\/}] R. E. Schwartz, {\it The Projective Heat Map\/},
A.M.S. Research Monograph, 2017 (to appear).
\newline
\newline
[{\bf SK\/}] E. B. Saff and A. B. J. Kuijlaars,
{\it Distributing many points on a Sphere\/}, 
Math. Intelligencer, Volume 19, Number 1, December 1997 pp 5-11
\newline
\newline
\noindent
[{\bf Th\/}] J. J. Thomson, {\it On the Structure of the Atom: an Investigation of the
Stability of the Periods of Oscillation of a number of Corpuscles arranged at equal intervals around the
Circumference of a Circle with Application of the results to the Theory of Atomic Structure\/}.
Philosophical magazine, Series 6, Volume 7, Number 39, pp 237-265, March 1904.
\newline
\newline
[{\bf T\/}] A. Tumanov, {\it Minimal Bi-Quadratic energy of $5$ particles on $2$-sphere\/}, Indiana Univ. Math Journal, {\bf 62\/} (2013) pp 1717-1731.
\newline
\newline
[{\bf W\/}] S. Wolfram, {\it The Mathematica Book\/}, 4th ed. Wolfram Media/Cambridge
University Press, Champaign/Cambridge (1999)
\newline
\newline
[{\bf Y\/}], V. A. Yudin, {\it Minimum potential energy of a point system of charges\/}
(Russian) Diskret. Mat. {\bf 4\/} (1992), 115-121, translation in Discrete Math Appl. {\bf 3\/} (1993) 75-81